\documentclass[11pt,letterpaper]{amsart}
\usepackage{amssymb, bm}
\usepackage{amsbsy}
\usepackage{amscd}
\usepackage{amsmath}
\usepackage{amsthm}
\usepackage{charter}
\usepackage[active]{srcltx}
\allowdisplaybreaks
\numberwithin{equation}{section}

\def\CC{{\mathbb C}}

\def\KZ{\mathbb{KZ}}
\def\LL{{\mathbb L}}

\def\QQ{{\mathbb Q}} 
\def\RR{{\mathbb R}} 
 
\def\ZZ{{\mathbb Z}} 

\def\Abold{\mathbf {A}}
\def\Ibold{\mathbf {I}}
\def\Pbold{\mathbf {P}}

\def\Sbold{\mathbf {S}}

\def\lambdabold{\bm{\lambda}}
\def\mbold{\mathbf {m}}
\def\zbold{\mathbf {z}}

\def\G{\Gamma}

\def\prim{{\rm prim}} 
\def\coprim{{\rm coprim}}

\newcommand{\p}{\partial}

\def\Bcal{{\mathcal B}}
\def\Ccal{{\mathcal C}}

\def\Ecal{{\mathcal E}} 
 
\def\Hcal{{\mathcal H}} 
\def\Ical{{\mathcal I}} 
\def\Jcal{{\mathcal J}} 

\def\Lcal{{\mathcal L}}
\def\Mcal{{\mathcal M}}
\def\Ocal{{\mathcal O}}

\def\Scal{{\mathfrak S}}

\def\Wcal{{\mathcal W}}

\def\half{{\tfrac{1}{2}}}

\def\zz{\mathbf {z}}

\def\glie{{\mathfrak{g}}}
\def\hlie{{\mathfrak{h}}}

\def\pt{{\scriptscriptstyle\bullet}}

\newcommand\ad{\operatorname{ad}}
\newcommand\aut{\operatorname{Aut}}
\newcommand\back{\operatorname{back}}

\newcommand\End{\operatorname{End}}

\newcommand\SL{\operatorname{SL}}

\newcommand\res{\operatorname{Res}}

\newcommand\slin{\operatorname{\mathfrak{sl}}}
\newcommand\orient{\operatorname{or}}

\newtheorem{theorem}{Theorem}[section]
\newtheorem{lemma}[theorem]{Lemma}
\newtheorem{proposition}[theorem]{Proposition}
\newtheorem{corollary}[theorem]{Corollary}

\newtheorem{lemmadef}[theorem]{Lemma-definition}

\newtheorem{conjecture}[theorem]{Conjecture}

\numberwithin{definition}{section}

\theoremstyle{remark}
\newtheorem{example}[theorem]{Example}
\newtheorem{remark}[theorem]{Remark}

\title[Polydifferentials and the KZ system]{A topological interpretation of the KZ system}

\author[Eduard Looijenga]{Eduard Looijenga}

\address{Mathematisch Instituut , Universiteit Utrecht,
PO Box 80.010, 3508 TA Utrecht, Nederland} 
\email{E.J.N.Looijenga@uu.nl}

\subjclass[2010]{Primary 32G34; Secondary 14D07}

\keywords{KZ system, polydifferentials, highest weight module}

\begin{document}

\begin{abstract}
We show that the KZ system has a topological interpretation in the sense that it may be understood as a variation of complex mixed Hodge structure whose successive pure weight quotients are polarized. This in a sense completes and elucidates work of Schechtman-Varchenko done in the early 1990's.  A central ingredient is a new realization of the irreducible highest weight representations of a Lie algebra of Kac-Moody type, namely on an algebra of rational polydifferentials on a countable product of Riemann spheres. We also obtain the kind of properties that  in the $\slin (2)$ case are due to Ramadas and are then known to imply  the unitarity of the WZW system in genus zero. 
\end{abstract}

\maketitle

\section*{Introduction}
We construct an identification of a given KZ-space with a space constructed out of cohomology with supports of a rank one local system. The latter space is topologically defined, but depends on the choice of $n$ distinct points on the affine line. By letting these points move, we get over the parameter space of such $n$-tuples a (Gau\ss-Manin) connection on the trivial bundle with fiber the  KZ-space and we show that this connection gets identified with the KZ-connection.
In a sense this completes earlier work of Schechtman-Varchenko, who constructed a flat map from the KZ-system to a  Gau\ss-Manin system for the ordinary cohomology (with no supports). Their map factors through the map constructed here, but 
they were not able to say much about how nontrivial that map might be in any given case (a priori it could be the zero map), whereas ours is always an isomorphism. We should perhaps emphasize that the topological interpretation tells us a great deal about the monodromy of the KZ-system that may be hard to obtain otherwise: the fundamental group of the base is a colored braid group and we can see it act on the cohomology of a rank one  local system over a space that is naturally attached to $n$ distinct points on the complex line. This space is like a partial  Deligne-Mumford compactification of the moduli space of genus zero curves that are endowed with $1+n+m$ punctures for some $m$ of which the first $1+n$ ones are fixed and the others freely move. We do not work this out here, but we feel that a closer examination of this representation is not only feasible, but also desirable.

Our approach differs from earlier work  in an other respect by putting at the center a bialgebra of rational polydifferentials on an iterated self-product of a Riemann sphere (these are like ordinary differential forms, except that this algebra, rather than obeying a Koszul rule, is plainly commutative). We show that if the number of factors is countably infinite, then such an algebra receives all the irreducible highest weight representations of all Lie algebras of Kac-Moody type, a fact that may have an interest in its own right. This makes for a smooth passage from the KZ system to the Gau\ss-Manin system and yields precious additional information. Thus the somewhat elusive KZ system is annexed by algebraic geometry, which  in the present setting means that the powerful tools of mixed Hodge theory become available to ferret out any finer structure. The strength of this approach manifests itself also in the algebraic domain, witness the new vanishing property Theorem \ref{thm:diagonaltensor} and the remarkable invariance property Theorem \ref{thm:slinvariances}.

We expect this set-up  to be particularly useful in establishing whether the WZW subsystem of the KZ-system has a flat unitary structure, as is conjectured by physicists and for which the evidence so far is limited to the case  $\slin (2)$ (due to Ramadas \cite{ramadas}, see also our version \cite{looijenga:sl2}). We characterize the WZW system in terms of the a simple vanishing property Corollary \ref{cor:wzwvanishing} that generalizes the one found by Ramadas. (Indeed, as is shown in 
\cite{looijenga:sl2} this result together with the vanishing theorem and the invariance property mentioned above  suffice  to obtain the unitarity for the  $\slin (2)$-case.)

While this paper may help to resolve one mystery, it creates another and that is the lack of a conceptual explanation for the appearance of polydifferentials. This is not the first place where these do occur in this context: they also do in a paper of Stoyanovski-Feigin \cite{sf} and in one by Beilinson-Drinfeld \cite{bd} (both are mainly concerned with the WZW system) and our feeling is that especially the first of these might  be linked to the present paper.

Although this article owes much to earlier work done in this area, the chosen approach almost demands it to be self-contained. For this reason we have given a complete proof that the vector bundle map between the vector bundles underlying the KZ system to the Gau\ss-Manin system is flat, despite the fact that the pattern of this proof is quite similar to (and indeed, inspired by) the one of Schechtman-Varchenko \cite{sv}.

Here is brief review of the separate sections. The first section is preparatory in nature:  we introduce here the basic algebra of polydifferentials that is at the center of this paper and derive some of its properties. This is used in Section 2 to show that this algebra contains the highest weight representations of a Lie algebra defined by a generalized Cartan matrix. We also make the transition of its  weight spaces to spaces of twisted logarithmic forms. In Section 3 we use the Casimir operator in order to define a Gau\ss-Manin connection and we show its compatibility  with the KZ connection: the main result in this direction, Theorem \ref{thm:gz=kz}, as well as the  version that yields a genuine topological interpretation of the KZ system, Theorem \ref{thm:bettermain}, in terms of cohomology with supports are stated here. 
As the underlying computations are somewhat delicate, we decided to write these out in full.  We state their core as an operator formula, but in order not to interrupt the flow of arguments, we relegated the proof to an appendix.  Our motivation was a desire to better understand the  WZW-systems  in genus zero as subsystems of the KZ-systems and Section 4 testifies to that. We characterize this system in terms of our polydifferentials and  
Conjecture \ref{conj:main} links them to a square integrability property. If true, then this  should realize  a WZW-system as a flat subbundle of polarized variation of Hodge structure which is of pure bidegree $(m,0)$ for some $m$  (so that this part is not really varying after all!) and would give the  flat unitary structure on this system  that physicists since long conjecture to exist.
\\

\textit{Conventions.}
For any finite sequence $I=(i_N,i_{N-1},\dots ,i_1)$  taken from an (index) set $\Ical$, we denote by $|I|:=N$ its length, by $\{ I\}:=\{ i_N,\dots ,i_1\}\subset \Ical$ the set of its terms  and by $I^*:=(i_1,i_{2},\dots ,i_N)$ the opposite sequence.
If  $(f_i)_{i\in\Ical}$ is an indexed collection of elements of an algebra  with unit, then $f_I$ denotes the noncommutative monomial  $f_{i_N}\dots f_{i_1}$ (understood to be $1$ if $I=\emptyset$).
Similarly, if $(f_i)_{i\in\Ical}$ is an indexed collection of elements of a  Lie algebra, then
$[f_I]$ stands for the iterated Lie bracket $\ad(f_{i_N})\ad(f_{k_{i-1}})\cdots \ad(f_{i_2})(f_{i_1})$ (read $f_i$  when $I=(i)$ and zero when $I=\emptyset$).  But if for a subset $X\subset\Ical$, $(c_x)_{x\in X}$ is a finite collection elements of an abelian group (usually complex numbers or elements of a vector space), then $c_X$ is sometimes used as an abbreviation for $\sum_{x\in X}c_x$ by $c_X$. If $t_i$ is a coordinate on a Riemann surface $C_i$ ($i\in\Ical)$, so that $(t_i)_{i\in\Ical}$ is a set of coordinates on the product $\prod_{i\in\Ical} C_i$, then any finite subset $X\subset \Ical$ defines a locus
where all $t_x$ with $x\in X$ are equal to each other; the resulting coordinate on that locus is the denoted by $t_X$ (the remaining coordinates  are $(t_i)_{i\in\Ical -X}$). Strictly speaking these conventions clash, but in practice confusion is unlikely.

The permutation group of a set $X$ is denoted by $\Scal (X)$.

\section{Some operators in an algebra of polydifferentials}\label{sect:polydiff}

We introduce an algebra of logarithmic polydifferentials which looks like a shuffle algebra and identify in it certain
operators of interest.

\subsection*{Polydifferentials} We consider polydifferentials on $C^\Ical$, where
$C$ is a nonsingular  complex curve and $\Ical$ is a finite set. A polydifferential of degree $d$ looks like a $d$-form on $C^\Ical$, but should not be confused with it as differentials coming from different factors commute rather than anti-commute. 

The definition is as follows. 
Let  $\Omega_C$ stand for the $\Ocal_C$-module of differentials and denote by $\Omega_C^\pt$ the sheaf of graded $\Ocal_C$-algebras $\Ocal_C\oplus\Omega_C$. For every $X\subset \Ical$, we have a natural, $\Scal (X)$-equivariant, identification of $\Ocal_{C^X}$-modules
\[
\Omega_C^{\boxtimes X}\cong\Omega^{|X|}_{C^X}\otimes\orient (X)
\]
where $\Omega_C^{\boxtimes X}$ denotes the exterior product $\otimes_{x\in X} \pi_x^*\Omega^1_C$ and 
$\orient(X)$ stands for the orientation module $\wedge^{|X|}\ZZ^X$. We define the \emph{sheaf of polydifferentials}  on $C^\Ical$ as the sheaf of graded $\Ocal_{C^{\Ical}}$-algebras $(\Omega^\pt_C)^{\boxtimes \Ical}$. So its degree $N$ part is 
$\bigoplus_{X\subset \Ical, |X|=N}\pi_X^*\Omega_C^{\boxtimes X}$,
where $\pi_X: C^{\Ical}\to C^X$ is the evident projection. A total order on $\Ical$ identifies the graded sheaf of polydifferentials with the graded sheaf of holomorphic differential forms on $C^{\Ical}$, but 
there is no natural analogue of the exterior derivative defined on  $(\Omega^\pt_C)^{\boxtimes \Ical}$.
If $\omega\in (\Omega^\pt_C)^{\boxtimes \Ical}$, we will write $\omega^X$ for its component in 
$\pi_X^*\Omega_C^{\boxtimes X}$, where $\pi_X: C^\Ical\to C^X$ is the evident projection. So $\omega=\sum_{X\subset \Ical} \omega^X$. We also write $\omega_X\in \pi_X^*(\Omega^\pt_C)^{\boxtimes X}$ for the sum $\sum_{Y\subset X}\omega^Y$.

We say that a meromorphic polydifferential on $C$  is \emph{logarithmic} if  locally it can be written as a product  
$(df_1/f_1)\dots (df_N/f_N)$ of logarithmic differentials , with $f_1,\dots ,f_N\in\Ocal_{C^\Ical}$.

Residue operators that involve a single factor (and a priori only those) have a meaning for polydifferentials.
Recall for any $p\in C$ is defined a residue operator 
$\res_{(t=p)}$ as a linear form on the space meromorphic differentials on $C$ at $p_i$.
This extends to polydifferentials: $\res_{(t_i=p)}$
maps the meromorphic polydifferentials on $\Ccal^\Ical$ to meromorphic polydifferentials on $\Ccal^{\Ical-\{i\}}$. But a residue along any other type of hypersurface must be treated with care, as it requires an ordering of the index set. For instance, if $(i,j)$ is an ordered pair, then the residue  along the diagonal $(t_i=t_j)$  where the $i$th and $j$th component are equal is to be understood as viewing $t_i$ as the residue variable and $t_j$ as a parameter, in other words, we take the residue in the $t_i$-direction and do this at the point with $t_i$-coordinate equal to $t_j$. The result is a polydifferential on the diagonal divisor $(t_i=t_j)$. In order to emphasize the asymmetric roles of $i$ and $j$, we shall denote this by $\res_{(t_i\to t_j)}$. This is clearly not the same thing as $\res_{(t_j\to t_i)}$.

\subsection*{A shuffle algebra of polydifferentials}

We now take $C$ to be the projective line  (which we will denote by $\Pbold$), we pick a point $\infty\in\Pbold$ and denote by $\Abold$ the affine line $\Pbold-\{\infty\}$. It is helpful to fix an affine coordinate on $\Abold$ (which depending on the context is denoted $t$ or $z$), although this is inessential  (in practice  it will only enter through differentials of the form $(t(p_1)-t(p_2))^{-1}dt$, where $p_1, p_1\in\Abold$ are distinct, which is indeed independent of $t$).

We consider relative polydifferentials on $(\Pbold^\Ical)_\Abold$, i.e., on  the projection 
$\Pbold^\Ical\times \Abold\to \Abold $. 
The coordinate of the base $\Abold$ is denoted $z$ and the affine coordinate of the $i$th factor $\Pbold$ in the product by $t_i$. 

We first consider the graded $\CC$-vector space $\Bcal_\Ical^\pt$ of the relative polydifferentials $(\Pbold^\Ical)_\Abold$ generated by the expressions of the form
\[
\zeta_I: z\mapsto \zeta_I(z):=\frac{dt_{i_N} dt_{i_{N-1}} \cdots dt_{i_1} }
{(t_{i_N}-t_{i_{N-1}})\cdots (t_{i_2}-t_{i_1})(t_{i_1}-z)}.
\]
where $I=(i_N,i_{N-1},\dots ,i_1)$ runs over the finite sequences in $\Ical$ (the degree of such a polydifferential is $N$). We stipulate that for $I=\emptyset$ we get $1$:
$\zeta_{\emptyset}=1$. Notice that $\zeta_I=0$ unless  the sequence $I$  is without repetition. In this paper $\Ical$ shall be  fixed (although we will later assume $\Ical$ to be countably infinite) and so we often write $\Bcal^\pt$ instead of $\Bcal_\Ical^\pt$.

We also put
\[
\omega_I=\frac{dt_{i_N} dt_{i_{N-1}} \cdots dt_{i_2} }
{(t_{i_N}-t_{i_{N-1}})\cdots (t_{i_2}-t_{i_1})},
\]
agreeing that when $I=(i)$, $\omega_{( i)}=\omega_i=1$ and that $\omega_{\emptyset}=0$.
So when $j\in\Ical-\{I \}$, then $\zeta_I(t_j)=\omega_{Ij}$ and when $I$ is nonempty and ends with $i$, then $\omega_I\zeta_{iJ}=\zeta_{IJ}$. 
It is clear that if $X\subset\Ical$, then $\Bcal_\Ical$ is closed under `taking the $X$-component'.

For a section $a$ of $\Pbold_\Abold$, the residue along $a$ is a map from the meromorphic differentials on $\Pbold_\Abold$ to $\CC(\Abold)$. Given an $x\in\Ical$, then this extends this in an obvious manner  as a map of degree $-1$ from to the meromorphic relative polydifferentials on $(\Pbold^\Ical)_\Abold$ to itself (strictly speaking it is the composite of a residue and the pull-back along the projection that suppresses the $x$th factor). For certain values of $a$, this preserves $\Bcal $:

\begin{lemmadef}\label{lemma:residues}
For $x\in\Ical$, we put 
\[
E_x:=-\res_{(t_x=\infty)}, \quad E'_x:=\res_{(t_x=z)}.
\]
Then for a sequence $I$ in $\Ical$ without repetition, we have $E_x(\zeta_I)=0$ resp.\ $E'_x(\zeta_I)=0$
 unless $I$ is of the form $xJ$ resp.\  $Jx$, in which case we get $\zeta_J$.

If $x,y\in \Ical$ are distinct,  then $\res_{(t_x\to t_y)}\zeta_I$ is zero unless
$I$ is of the form $I''xyI'$ or $I''yxI'$, in which case we get $\zeta_{I''yI'}$ resp.\ 
$-\zeta_{I''yI'}$.
\end{lemmadef}

The proof is straightforward.
This helps us (among other things) to give an intrinsic characterization of 
$\Bcal^\pt$. 
Given a subset $X\subset\Ical$, denote by $D_X$ the reduced effective divisor on $\Pbold^X\times\Pbold$ that is the sum of the loci $(t_i=\infty)$, $(t_i=z)$ and the diagonals $(t_i=t_j)$ ($i,j\in X$ distinct).

\begin{corollary}\label{cor:intrinsic}
The polydifferentials $\zeta_I$ are linearly independent over $\CC [z]$ and (hence) constitute a $\CC$-basis of $\Bcal^\pt$. Moreover, we have a decomposition $\Bcal^N=\oplus_{X\subset\Ical, |X|=N} \Bcal_X^N$ and when $X\subset \Ical$ is finite,
then the natural map of $\CC[z]$-modules
\[
\CC[z]\otimes_\CC \Bcal_X^{|X|}\to H^0(\Pbold^X\times\Abold,\Omega^{|X|}_{\Pbold^X \times\Abold/\Abold}(\log D_X))\otimes\orient (X)
\]
is an isomorphism that  maps $\Bcal^{|X|}$ onto the space of sections that vanish at $z=\infty$ and  restricts for any $z\in\Abold$ to an isomorphism
\[
\Bcal_X^{|X|}\to H^0(\Pbold^X,\Omega^{|X|}_{\Pbold^X}(\log D_X(z)))\otimes\orient (X).
\]
\end{corollary}
\begin{proof}
The first assertion follows from the observation that for  a sequence $I=(i_N,\dots ,i_1)$, the iterated residue $E'_{i_N}\cdots E'_{i_1}$ takes on $\zeta_J$ the value  $\delta_{I,J}$ when $|J|\le |I|$. 

It is clear that $\Bcal_\Ical^\pt=\oplus_{X\subset\Ical} \Bcal_X^{|X|}$ and it is easily checked  that  $\Bcal_X^N$ is a subspace of 
$H^0(\Pbold^X\times\Pbold,\big(\Omega^{\boxtimes X}_{\Pbold}\big)_\Pbold(D_X-(z=\infty)))$. To see that this inclusion is an equality, we observe that if an element $\omega$  of the last space is annihilated by all the iterated residue operators of the type above, then $\omega$ must be  independent of $z$. As $\omega$ is zero for $z=\infty$, it follows that $\omega=0$.
The last assertion is proved similarly.
\end{proof}

Here is another interesting consequence of Lemma \ref{lemma:residues}. 
Given a sequence $I=(i_N,\dots ,i_1)$ in $\Ical$ without repetition of length $N\ge 2$, then put  
\[
\res_I:=\res_{(t_{i_2}\to t_{i_1})}\cdots \res_{(t_{i_N}\to t_{i_{N-1}})}.
\]
We regard this operator as taking values in the polydifferentials on the diagonal locus defined by $t_{i_1}=t_{i_2}=\cdots =t_{i_N}$.
We recall that $t_{\{I\}}$ denotes the restriction of any $t_{i_k}$ to this locus.

\begin{lemma}\label{lemma:flagresidue}
Let $I=(i_N,\dots ,i_1)$ be a sequence without repetition in $\Ical$. Then the following identities hold for $a\in\{ z,\infty\}$ in  $\Bcal $:
\begin{align*}
[[[\cdots [E_{i_N},E_{i_{N-1}}],\cdots ] ,E_{i_2}],E_{i_1}] &=-\res_{(t_{\{ I\}}=\infty)}\res_I,\\
[[[\cdots [E'_{i_N},E'_{i_{N-1}}],\cdots ] ,E'_{i_2}],E'_{i_1}] &=\res_{(t_{\{ I\}}=z)}\res_I.
\end{align*}
\end{lemma}
\begin{proof} 
The last clause of Lemma \ref{lemma:residues} shows that if $\res_I\zeta_{K}\not= 0$, then some permutation $\sigma (I)$ of $I$ must appear as an uninterrupted subsequence of $K$: $K=K''JK'$, and  have the following property: $i_N$ is adjacent to $i_{N-1}$ and removing $i_N$ makes $i_{N-1}$ adjacent to $i_{N-2}$ and so on. In other words,  there is a proper subsequence $I'< I$ such that if $I''$ is the residual subsequence (obtained by removing the terms in $I'$), then $\sigma (I)=(I')^*I''$;  let us call this a \emph{back-forward} permutation of $I$ and  denote by $\back (\sigma)$ the length of $I'$. If we also require that 
$\res_{(t_{\{ I\}}=z)}\res_I (\zeta_{K})\not=0$, then $K'=\emptyset$ and its value is then by
Lemma \ref{lemma:residues} equal to  $(-1)^{\back (\sigma)}\zeta_{K''}$. If on the other hand 
$-\res_{(t_{\{ I\}}=\infty)}\res_I (\zeta_{K})$ is nonzero, then $K''=\emptyset$ and its value is then by
Lemma \ref{lemma:residues} equal to $(-1)^{\back (\sigma)}\zeta_{K'}$. 

Expanding the iterated bracket of residue operators in a straightforward fashion yields
\begin{multline*}
[[[\cdots [E'_{i_N},E'_{i_{N-1}}],\cdots ] ,E'_{i_2}],E'_{i_1}] =\\
=\sum_\sigma (-1)^{\back (\sigma)}E'_{\sigma(i_N)}E'_{\sigma(i_{N-1})}\cdots E'_{\sigma(i_1)},
\end{multline*}
where the sum is over all back-forward  permutations $\sigma$ of $I$. Its value on $\zeta_K$ is nonzero only if $K$ is of the form $K''\sigma (I)$  and is then equal to $\zeta_{K''}$. If we replace  $E'_{i_k}$
by $E_{i_k}$ and use the antisymmetry of the bracket: then we find expansion above is also equal to
\begin{multline*}
[[[\cdots [E_{i_N},E_{i_{N-1}}],\cdots ] ,E_{i_2}],E_{i_1}] =\\=
\sum_\sigma (-1)^{\back (\sigma)}E_{\sigma (i_1)}E_{\sigma(i_2)}\cdots E_{(t_{\sigma (i_N)}}.
\end{multline*}
Its value on is $\zeta_K$ is nonzero only if $K$ is of the form $\sigma (I)K'$ and is then equal to $(-1)^{\back(\sigma)}\zeta_{K'}$.
\end{proof}

A subset $X\subset\Ical$ induces via pull-back an injective map from rational polydifferentials on $C^X$ to rational polydifferentials on $C^\Ical$.  We now allow $\Ical$ to be infinite, but countable and define the space of polydifferentials on $C^\Ical$ as the injective limit of the polydifferentials on $C^X$, where $X$ runs over the finite subsets of $\Ical$. 
We denote by $\hat\Bcal^m$ the space of  (possibly infinite) sums of these relative polydifferentials of degree $m$:
\[
\hat\Bcal^m=\prod_{X\subset\Ical, |X|=m} \Bcal_X^m
\] 
and put $\hat\Bcal^{\pt}:= \oplus_{m=0}^\infty\hat\Bcal^m$, but often omit the subscript $\Ical$ as that has been fixed. 

If $I$ and $J$ are sequences such that their juxtaposition $IJ$ is without repetition, then we denote by 
$I\star J$ the collection of \emph{shuffles} of $I$ and $J$, i.e., the set of sequences that are permutations of $IJ$ and in which
$I$ and $J$ appear as subsequences. This \emph{shuffle product} is associative: we have that 
$(I\star J)\star K=
I\star (J\star K)$ if $IJK$ is without repetition.

\begin{lemma}[Shuffle rules]\label{lemma:shuffle}
The graded vector space $\hat\Bcal$ is closed under product and is as a $\CC$-algebra isomorphic to a (completed) shuffle algebra: if $I$ and $J$ are finite sequences in $\Ical$, then  $\zeta_I\zeta_J=\sum_{K\in I\star J}\zeta_{K}$ (we get zero unless $IJ$ is without  repetition). Together with the coproduct
\[
\delta :\hat\Bcal\to \hat\Bcal\otimes \hat\Bcal,\quad \delta (\zeta_I)=\sum_{I=I''I'} \zeta_{I''}\otimes \zeta_{I'},
\]
this makes $\hat\Bcal$ a commutative bialgebra over $\CC$ (with the projection on the degree zero summand as counit).

More generally, if $I$ and $J$ are finite nonempty sequences in $\Ical$,  $k\in \{J\}$ and 
 $J$ is written $J_{<k}J_{\ge k}$, then 
$\omega_{Ik}\zeta_J=\sum_{K\in I\star J_{<k}}\zeta_{KJ_{\ge k}}$.
\end{lemma}
\begin{proof} 
An induction argument shows that it is enough to verify these statements insofar they do not regard the Hopf property when $I$ is a singleton. That case easily follows from repeated use from the simple identity 
\[
\frac{1}{v-u}-\frac{1}{w-u}=\frac{w-v}{(v-u)(w-u)}.\qedhere
\]
The proof that $\delta$ is a coproduct that is compatible with the shuffle product is straightforward.
\end{proof}

So $\Bcal$ can be identified with a shuffle algebra over $\CC[t]/(t^2)$ on a set of generators indexed by $\Ical$. This makes us wonder whether there is a relation with iterated integrals. From Corollary \ref{cor:intrinsic} and 
Lemma \ref{lemma:shuffle} we deduce:

\begin{corollary}\label{cor:polynomial}
If $\Ical$ is a countably infinite set, then $\hat\Bcal^{\Scal (\Ical)}$ is a polynomial algebra  
with primitive generator $\zeta:=\sum_{i\in\Ical} \zeta_{i}$ and we have
\[
\frac{\zeta^N}{N!}=\sum_{\{X\subset \Ical, |X|=N\}} \prod_{x\in X}\frac{dt_x}{t_x-z}.
\]
\end{corollary}

\subsection*{Two-parameter identities}
The following two lemmas assert identities in an algebra of meromorphic 
polydifferentials on $\Pbold$ depending on two complex  variables.
They will be needed later (beginning with the construction of the Gau\ss-Manin connection in
Lemma \ref{lemma:gm}), but it is convenient to state and prove them now. 

Recall that for a sequence $I$, $I^*$ denotes the opposite sequence.

\begin{lemma}\label{lemma:mixedshuffle}
Let $I$ be a nonempty sequence in $\Ical$. If $b(I)$ denotes its first element, then
\[
\frac{z-w}{t_{b(I)}-w}\zeta_I(z)
=\sum_{I=I_2I_1} (-1)^{|I_2|}\zeta_{I_1}(z)\cdot\zeta_{I_2^*}(w),
\]
or equivalently, the operator $\sum_{J} (-1)^{|J|}\zeta_{J} (w)E_J$ (with the sum taken over all finite sequences $J$ in $\Ical$) sends  $\zeta_I$ to $\frac{z-w}{t_{b(I)}-w}\zeta_I$.

More generally, if $i$ appears in $I$ so that we can write $I=I''iI'$, then 
\[
\frac{z-w}{t_i-w}\zeta_I(z)
=\zeta_{I}(z)
-\sum_{I'=I_2I_1} (-1)^{|I_2|}\zeta_{I_1}(z)\cdot
\omega_{I''i}\zeta_{I_2^*i}(w),
\]
and if $j\in\Ical$,  $j\not= i$, then
\[
\frac{z-t_j}{t_i-t_j}\zeta_I(z)
=\zeta_{I}(z)
-\sum_{I'=I_2I_1} (-1)^{|I_2|}\zeta_{I_1}(z)\cdot
\omega_{I''i}\omega_{I_2^*ij}.\\
\]
\end{lemma}
\begin{proof}
The third identity is not really different from the second (take $t_j=w$)  and  the second follows from the first applied to $iI'$  (and multiply it with $\omega_{I''i}$). So we concentrate on the first identity and prove it with induction on $|I|$. For $I=\emptyset$ there is nothing to show and for $I=(i)$, the lemma states that 
\[
\frac{z-w}{t_i-w}\zeta_i(z)=\zeta_{i}(z)-\zeta_i(w),
\]
which is a simple consequence of the identity
$\frac{z-w}{t_i-w}=1-\frac{t_i-z}{t_i-w}$. Now assume  $I$ has length $>1$ and write $I=iI'$ and $I'=jJ'$.  
So $\zeta_I(z)=\zeta_{ij}\zeta_{jJ}(z)=
\frac{dt_i}{t_i-t_j}\zeta_{I'}(z)$. Since
\[
\frac{z-w}{(t_{i}-w)(t_{i}-t_{j})}=\frac{z-w}{t_{j}-w}\Big( \frac{1}{t_{i}-t_{j}}-\frac{1}{t_{i}-w}\Big),
\]
we have
\begin{multline*}
\frac{z-w}{t_{i}-w}\zeta_I(z)=\frac{z-w}{(t_{i}-w)(t_{i}-t_{j})}dt_i
\zeta_{I'}(z)\\
=\frac{z-w}{t_j-w}\Big( \frac{dt_{i}}{t_{i}-t_{j}}-\frac{dt_{i}}{t_{i}-w}\Big)
\zeta_{I'}(z)\\
=\frac{z-w}{t_{j}-w}\zeta_I(z)-\frac{z-w}{t_{j}-w}\zeta_{I'}(z)\cdot\zeta_{i}(w),
\end{multline*} 
which after invoking the induction hypothesis becomes
\begin{multline*}
\zeta_I(z)-\sum_{J'=J_2J_1} (-1)^{|J_2|}\zeta_{J_1}(z)\cdot
\omega_{J_2^*j}\zeta_{ij}(w)\\
-\Big(\zeta_I'(z)-\sum_{J'=J_2J_1}(-1)^{|J_2|}\zeta_{J_1}(z)\cdot\zeta_{J_2^*j}(w)\Big)\zeta_{i}(w)
=\zeta_I(z)+\\+\sum_{J'=J_2J_1} (-1)^{|J_2|}\zeta_{J_1}(z)\Big(-\omega_{J_2^*j}\zeta_{ij}(w)
+\zeta_{J_2^*j}(w)\zeta_{i}(w)\Big)-\zeta_{I'}(z)\cdot\zeta_{i}(w). 
\end{multline*} 
According to the shuffle rules \ref{lemma:shuffle}, we have
$-\omega_{J_2^*j}\zeta_{ij}(w)+\zeta_{J_2^*j}(w)\zeta_{i}(w)=
\zeta_{J_2^*ji}(w)= \zeta_{(jJ_2)^*i}(w)$, so that we get
\[
\zeta_I(z)-\sum_{J'=J_2J_1} (-1)^{|jJ_2|}\zeta_{J_1}(z)\zeta_{(jJ_2)^*i}(w)
-\zeta_{I'}(z)\cdot\zeta_{i}(w),
\]
which indeed equals $\zeta_I(z)-\sum_{I'=I_2I_1} (-1)^{|I_2|}\zeta_{I_1}(z)\zeta_{I_2^*i}(w)$.
\end{proof}

\begin{remark}
Notice that the expression $\omega_{I_2^*i}\zeta_{I''i}(w)$ in the right hand side 
of the preceding lemma is,  according to Lemma \ref{lemma:shuffle},  equal to the sum $\sum_L \zeta_{Li}(w)$, where $L$ runs over the shuffles of $I_2^*$ and $I''$.
\end{remark}

The following lemma generalizes the preceding one.

\begin{lemma}\label{lemma:mixedshuffle2}
If  $I:=I''iI'$ and $J:=J''jJ'$ are sequences as above, then
\begin{multline*}
\frac{z-w}{t_{i}-t_{j}}\zeta_I(z)\cdot\zeta_J(w)=
\zeta_I(z)\cdot\zeta_J(w)\\
-\sum_{I'=I_2I_1} (-1)^{|I_2|}\zeta_{I_1}(z)\cdot
\omega_{I_2^*i}\omega_{I''ij}\zeta_J(w)\\
-\sum_{J'=J_2J_1} (-1)^{|J_2|}
\omega_{J_2^*j}\omega_{J''ji}\zeta_I(z)\cdot\zeta_{J_1}(w).
\end{multline*}
In particular, it is a linear combination of terms of the form $\zeta_K(z)\cdot\zeta_L(w)$.
\end{lemma}
\begin{proof}
The second identity of Lemma \ref{lemma:mixedshuffle} gives after 
multiplication by $\zeta_J(w)$
\[
\frac{z-t_{j}}{t_{i}-t_{j}}\zeta_I(z)\cdot\zeta_J(w)
=\zeta_I(z)\cdot\zeta_J(w)
-\sum_{I'=I_2I_1} (-1)^{|I_2|}\zeta_{I_1}(z)\cdot
\omega_{I_2^*i}\cdot\omega_{I''ij}\zeta_J(w).
\]
Likewise we find 
\[
\frac{w-t_i}{t_{j}-t_i}\zeta_I(z)\cdot\zeta_J(w)=\zeta_I(z)\cdot\zeta_J(w)
-\sum_{J'=J_2J_1} (-1)^{|J_2|}
\omega_{J_2^*j}\omega_{J''ji}\zeta_I(z)\cdot\zeta_{J_1}(w).
\]
The lemma then follows by adding these two identities and using that
\[
\frac{z-t_{j}}{t_i-t_{j}} +\frac{w-t_i}{t_{j}-t_i}=\frac{z-w}{t_i-t_{j}}+1.\qedhere
\]
\end{proof}

\subsection*{The $\Phi$-operators}
In the space of rational relative polydifferentials on $(\Pbold^\Ical)_\Abold$, we regard  $dt_i$  ($i\in I$) not just as an element, but also as the multiplication operator in this space. Its adjoint acts in the $i$th tensor factor only and sends $dt_i$ to $1$ and $1$ to zero, hence is the contraction operator $\iota_{\p/\p t_i}$.

Let now be given complex numbers $(p_i)_{i\in \Ical}$ and $(c_{i,j})_{i,j\in\Ical, i\not=j}$. We define for $i\in\Ical$ an operator $\Phi_i$ in the space of relative meromorphic polydifferentials by 
\[
\Phi_i :=\frac{p_idt_i}{t_i-z} -\sum_{j\not=i} c_{i,j}\frac{dt_idt_j}{t_i-t_j}\iota_{\p/\p t_j},
\]
So for a finite subset $X\subset\Ical$, we have
\[
\Phi_i (dt_X)=\Big(\frac{p_i}{t_i-z}-\sum_{x\in X}\frac{c_{i,x}}{t_i-t_x}\Big)dt_idt_X,
\]
where it is understood that the right hand side is zero when $i\in X$.

\begin{lemma}\label{lemma:betaop}
This operator preserves $\Bcal $ (hence also  the completion $\hat\Bcal $), for
\[
\Phi_x (\zeta_I)=\sum_{I=I''I'}(p_x-c_{x,I'})\zeta_{I''xI'}.
\]
(Observe that the right hand side vanishes if $x\in I$ and that the  term indexed by $(I'',I')=(I,\emptyset)$ reduces to $p_x\zeta_{Ix}$.) Furthermore, for any $o\in \Ical$,
\[
[\Phi_x ,\omega_{Io}]=\sum_{I=I''I'} -c_{x,I'}\omega_{I''xI'o}.
\]
\end{lemma}
\begin{proof}
We have $\Phi_x (\zeta_I)=p_x\zeta_x\zeta_I-\sum_{i\in \{I\}} c_{x,i}\omega_{x,i}\zeta_I$, and the latter is by Lemma \ref{lemma:shuffle} equal to 
\[
p_x \sum_{I=I''I'}\zeta_{I''xI'}-\sum_{I=I''I'}\sum_{i\in \{ I'\}} c_{x,i}\zeta_{I''xI'}=
\sum_{I=I''I'}(p_x-c_{x,I'})\zeta_{I''xI'}.
\]
The second identity follows from this.
\end{proof}

\subsection*{Products and powers} Notice that  
\[
[\Phi_i ,dt_X]=\sum_{x\in X} \frac{-c_{i,x}dt_i}{t_i-t_x}dt_X.
\]
We also observe that for a sequence $I=(i_N,\dots ,i_1)$ in $\Ical$
\begin{multline*}
\Phi_I( dt_X):=\Phi_{i_N}\Phi_{i_{N-1}}\cdots \Phi_{i_1} (dt_X)\\
 =\prod _{I=I''iI'} dt_{i}\Big(\frac{p_{i}}{t_{i}-z}-\sum_{x\in X} \frac{c_{i,x}}{t_{i}-t_x} -\sum_{j\in I'} \frac{c_{i,j}}{t_{i}-t_{j}}\Big).dt_X\\
= \prod _{I=I''iI'}\frac{dt_{i}}{t_{i}-z}\Big(p_{i}
 - \sum_{x\in X} c_{i,x}\frac{t_{i}-z}{t_{i}-t_x} - \sum_{j\in \{ I'\}} c_{i,j}\frac{t_{i}-z}{t_{i}-t_{j}}\Big).dt_X.
\end{multline*}

\begin{lemma}\label{lemma:sumproductidentity}
In $\CC[a,b]$ we have
\[
\sum_{\sigma\in\Scal_N}\prod_{k=1}^N \Big(b+a\sum_{l<k}\frac{t_{\sigma(k)}-z}{ t_{\sigma(k)}-t_{\sigma(l)}}\Big)=
N! b(b+\half a)\cdots (b+\half (N-1)a).
\]
\end{lemma}
\begin{proof}
The lemma is evidently true for $N=1$. We continue with induction on $N$. Let us write 
$c_N$ for $N! b(b+a)\cdots (b+(N-1)a)$. 
Then the left hand side is equal to
\begin{multline*}
\sum_{k=1}^N\sum_{\substack{\sigma\in\Scal_N\\ \sigma (N)=k}}\prod_{k=1}^N \Big(b+a\sum_{l<k}\frac{t_{\sigma(k)}-z}{ t_{\sigma(k)}-t_{\sigma(l)}}\Big)=\\
=\sum_{k=1}^N
\Big(b+a\sum_{l\not=k}\frac{t_k-z}{ t_k-t_l}\Big)\sum_{\substack{\sigma\in\Scal_N\\ \sigma (N)=k}}
\prod_{k=1}^{N-1} \Big(b+a\sum_{l<k}\frac{t_{\sigma(k)}-z}{ t_{\sigma(k)}-t_{\sigma(l)}}\Big)=\\
=\sum_{k=1}^N \Big(b+a\sum_{l\not=k}\frac{t_k-z}{ t_k-t_l}\Big)c_{N-1}=
\Big(Nb+a\sum_{l<k}\big(\frac{t_k-z}{ t_k-t_l}+\frac{t_l-z}{ t_l-t_k}\big) \Big)c_{N-1}=\\
=(Nb+a\sum_{l<k} 1)c_{N-1}
=Nc_{N-1}\big(b +\half (N-1)a\big)=c_N.\qedhere
\end{multline*}
\end{proof}

Of special interest is the case when $b$ is a positive integer $m$ and $a=-2$. Then we find
that the right hand side of Lemma \ref{lemma:sumproductidentity} is zero for $N>m$ and equals $\frac{N! m!}{(m-N)!}$ otherwise. The lemma above and the discussion preceding it show:

\begin{corollary}\label{cor:power}
Let $\Jcal\subset \Ical$  be such that for all $i,j\in\Jcal$ we have $c_{i,j}=2$ and $p_i=m$ for a fixed nonnegative integer $m$.  If $\Phi:=\sum_{i\in\Jcal} \Phi_i$, then
\[
\frac{\Phi^N}{N!}(1)= m(m-1)\cdots (m+1-N)\sum_{\{K\subset \Jcal: |K|=N\}}
\prod _{k\in \{K\}} \frac{dt_k}{t_k-z},
\]
where $K$ runs over all the $N$-element subsets of $\Ical'$. In particular, this expression vanishes for $N=m+1$.
\end{corollary}

\subsection*{Commutators}
A straightforward  check shows:

\begin{lemma}\label{lemma:simplebracket}
If $x\not=y$, then
$[\Phi_x,\Phi_y]= -c_{x,y}\omega_{x,y}\Phi_y+c_{y,x}\omega_{y,x}\Phi_x$.
\end{lemma}

We use this to prove:

\begin{lemma}\label{lemma:bracketidentity}
Let $I$ be a nonempty sequence in $\Ical$ without repetition and let $o\in \Ical$  not occur in $I$. Assume that for some complex numbers $a, c,\check{c}$ we have $c_{i,j}=a$ for all $i\not=j$ in $\{ I\}$ and 
$c_{i,o}=c$, $c_{o,i}=\check{c}$ for all $i\in \{ I\}$. 
If $\ell (I)$ denotes the last element of $I$, then we have the operator identity 
\[
[\Phi_{Io}]
=\prod_{\ell (I)\not=k\in \{ I\}}\Big(\frac{c}{t_o-t_k} +\sum_{l\in \{ I_{>k}\}}\frac{a}{t_l-t_k}\Big)\cdot
\frac{cdt_I\Phi_o+\check{c}\sum_{i\in \{I\}}dt_{I_i}\Phi_i}{t_o-t_{\ell (I)}},
\]
where $I_i$ is obtained from $I$ by omitting $i$.
\end{lemma}
\begin{proof}
Lemma \ref{lemma:simplebracket} gives the asserted identity for $I=(x)$: we have 
\[
[\Phi_x,\Phi_o]=\frac{1}{t_o-t_x}\big(cdt_x\Phi_o+\check{c}dt_{o}\Phi_{i}\big).
\]
This verifies the lemma when $I$ is a singleton. We proceed with induction on the length of $I$. So if we write $I=xJ$, then 
\begin{multline*}
[\Phi_{Io}]=[\Phi_x,[\Phi_{Jo}]]=\\
\prod_{\ell (I)\not=k\in J}
\Big( \frac{c}{t_o-t_k} +\sum_{l\in \{J_{>k}\}}\frac{a}{t_l-t_k}\Big)\cdot
\frac{c[\Phi_x,dt_J\Phi_o] +\check{c}\sum_{j\in J} [\Phi_x,dt_{J_jo}\Phi_h]}{t_o-t_{\ell (I)}}.
\end{multline*}
So the induction step amounts to showing that the numerator of the last factor,
$c[\Phi_x,dt_J\Phi_o] +\check{c}\sum_{j\in J} [\Phi_x,dt_{J_jo}\Phi_h]$, is equal to 
\[
\Big(\frac{c}{t_o-t_x}+\sum_{l\in J}\frac{a}{t_l-t_x}\Big)
\Big(cdt_{I}\Phi_o + \check{c}\sum_{i\in \{J\}}dt_{I_io}\Phi_i+\check{c}dt_{Jo}\Phi_x\Big).
\]
To this end we expand the brackets in the left hand side using the already verified case $N=1$:
\begin{multline*}
[\Phi_x,dt_J\Phi_o]=[\Phi_x,dt_J]\Phi_o+dt_J[\Phi_x,\Phi_o]=\\
=\sum_{l\in \{ J\}}a\frac{dt_x}{t_l-t_x}dt_J\Phi_o+\frac{dt_J}{t_o-t_x}\big(cdt_x\Phi_o+
\check{c}dt_o\Phi_x\big)\\
=\Big(\frac{c}{t_o-t_x} +\sum_{l\in \{ J\}}\frac{a}{t_l-t_x}\Big)dt_I\Phi_o
+\check{c}\frac{dt_{Jo}}{t_o-t_x}\Phi_x
\end{multline*}
and for $i\in J$, 
\begin{multline*}
[\Phi_x, dt_{J_io}\Phi_i]=[\Phi_x,dt_{J_io}] \Phi_i+dt_{J_io}[\Phi_x,\Phi_i]=\\
=\sum_{l\in \{J_i\}} a\frac{dt_{xJ_io}}{t_l-t_x}\Phi_i+ 
c\frac{dt_{xJ_io}}{t_o-t_x}\Phi_i+
a\frac{dt_{J_io}}{t_i-t_x}(dt_x\Phi_i+dt_i\Phi_x)=\\
=\Big(\frac{c}{t_o-t_x}+\sum_{l\in \{ J\}}\frac{a}{t_l-t_x}\Big)dt_{I_io}\Phi_i 
+a\frac{dt_{Jo}}{t_i-t_x}\Phi_x.
\end{multline*}
If we substitute these identities in the linear combination $c[\Phi_x,dt_J\Phi_o] +\check{c}\sum_{j\in J} [\Phi_x,dt_{J_jo}\Phi_h]$ we get the desired expression.
\end{proof}

For any $\Jcal\subset \Ical$, we put  $\tau_\Jcal:= \sum_{x\in \Jcal} dt_x$, so that
\[
\frac{\tau_\Jcal^N}{N!}=\sum_{X\subset \Jcal, |X|=N} dt_X\quad \text{(read  $1$ if  $N=0$)}.
\]

\begin{corollary}\label{cor:serretype}
Let $\Jcal\subset \Ical$, $o\in\Ical-\Jcal$, $c,\check{c}\in \CC$ be such that for all $x,y\in \Jcal$,
$c_{x,o}=c$, $c_{o,x}=\check{c}$ and $c_{x,y}=2$. If $\Phi:=\sum_{x\in \Jcal} \Phi_x$, then
\[
\frac{(\ad \Phi)^N}{N!}(\Phi_o)=(-c-1)\cdots (-c-(N-1))
\big(c\frac{\tau_{\Jcal}^N}{N!}\Phi_o  +\check{c}dt_o\frac{\tau_{\Jcal}^{N-1}}{(N-1)!}\Phi \big).
\]
In particular, the left hand side is zero for $N=-c+1$.
\end{corollary}
\begin{proof}
First observe that $(\ad \Phi)^N(\Phi_o)=\sum_{I} [\Phi_{Io}]$, 
where $I$ runs over  the sequences in $\Jcal$ of length $N$ without repetition.  Now apply Lemma's \ref{lemma:sumproductidentity}
and \ref{lemma:bracketidentity} to each $\Scal_N$-orbit of this index set.
\end{proof}

Lemma \ref{lemma:residues} yields:

\begin{lemma}\label{lemma:efcommutator}
Given $x,y\in \Ical$, then $[E_x,\Phi_y](\zeta_I)=0$ unless $x= y\notin I$, in which case it multiplies  $\zeta_I$ by the scalar  $p_x-c_{x,I}$.
\end{lemma}
\begin{proof} 
First we notice that  $\Phi_yE_x(\zeta_I)$ vanishes unless $I$ is of the form
$xJ$ and $x\notin J$: we  then get $\sum_{J=J''J'}(p_y-c_{y,J'})\zeta_{J''x J'}$. 

On the other hand $E_x\Phi_y(\zeta_I)$ vanishes unless $y\notin I$ and 
either $I$ has the form $xJ$ or $x=y$. In the first case, we get $\sum_{J=J''J'}(p_y-c_{y,J'})\zeta_{J''yJ'}=\Phi_yE_x\zeta_I$ and so $\zeta_I$ is killed by $[E_x,\Phi_y]$. In the second case, we get $(p_x-c_{x,I})\zeta_I$ and we note that then also $\Phi_yE_x(\zeta_I)=0$.
\end{proof}

\section{Highest weight representations in spaces of polydifferentials}\label{sect:reps}
In this section we show among other things that the highest weight representation of Lie algebras of Kac-Moody type are naturally realized in an algebra of logarithmic polydifferentials.

\subsection*{Kac-Moody Lie algebras} Let $(c_{k,l})_{k,l=1}^r$ be a generalized Cartan matrix, i.e.,  $c_{k,k}=2$, and  for $k\not= l$,  $c_{k,l}$ is a nonpositive integer which is zero if and only if $c_{l,k}$ is zero. Attached to this matrix is the  Lie algebra  defined by the following presentation: it has  generators $\tilde e_1,\dots ,\tilde e_r,\tilde f_1,\dots ,\tilde f_r$  subject to the relations $[\tilde e_k,\tilde f_l]=0$ for $k\not=l$ and if we put $\check{\alpha}_k:=[\tilde e_k,\tilde f_k]$, then
\[
[\check{\alpha}_k,\tilde e_l]=c_{k,l} \tilde e_l,\quad  [\check{\alpha}_k,\tilde f_l]=-c_{k,l} \tilde f_l,\quad   [\check{\alpha}_k,\check{\alpha}_l]=0.
\]
We define the Lie algebra $\glie$ as a quotient of  this Lie algebra by also imposing  the \emph{Serre relations} by setting for $k\not=l$,  $\ad (\tilde e_k)^{1-c_{k,l}}\tilde e_l$ and $\ad (\tilde f_k)^{1-c_{k,l}}\tilde f_l$ equal to zero. We denote  the linear span of the $\check{\alpha}_k$'s by $\hlie$. (We obtain the 
(Kac-Moody) Lie algebra as defined in \cite{kac} as  the quotient of $\glie$ by the maximal ideal of that has zero intersection with  $\hlie$, but as is shown in \emph{op.\ cit.}, we have equality in case the generalized Cartan matrix is symmetrizable, a condition that is always fulfilled in the cases of interest.) We denote by
$\tilde\glie$ the intermediate Lie algebra defined by imposing the latter half of these relations only:
so we let $\ad (\tilde f_k)^{1-c_{k,l}}\tilde f_l=0$. The images of $\tilde e_k$, $\tilde f_k$ and $\check{\alpha}_k$ in $\tilde\glie$ are denoted by the same symbol (so that is a slight abuse of notation), but in $\glie$ the first two lose their tilde's. The linear span of the $\check{\alpha}_k$'s, which we  shall denote by $\hlie$,  will be regarded as a  subalgebra of both $\tilde\glie$ and $\glie$. It is a Cartan subalgebra of either. Notice that  the simple root $\alpha_l\in\hlie^*$, characterized by $[h,e_l]=\alpha_l(h)e_l$, takes on
$\check{\alpha}_k$ the value $c_{k,l}$. We denote by $\tilde\glie_+\subset\tilde\glie$ the subalgebra generated by the
$\tilde e_k$'s and by $\tilde\glie_-\subset\tilde\glie$ the subalgebra generated by the $\tilde f_k$'s.

Let $\tilde V$ be a representation of $\tilde g$ on which $\hlie$ acts semisimply (and hence is graded by weights). The \emph{primitive} part $\tilde V^\prim\subset \tilde V$ is by definition the set of vectors killed by $\tilde\glie_+$ (the biggest subspace on which $\tilde\glie_+$ acts trivially), whereas the 
\emph{coprimitive} part of $V$ is the quotient $\tilde V_\coprim:=\tilde V/\tilde\glie_-V$ of $\tilde V$ (the smallest on which the $\tilde\glie_-$ acts trivially). Notice that both inherit a semisimple $\hlie$-action. The following lemma collects a few simple, but useful facts about the representation $\tilde V$.

\begin{lemma}\label{lemma:primg}
The $\tilde\glie$-submodule of $\tilde V$ generated by its primitive part $\tilde V^\prim$  is in fact a $\glie$-submodule. If  $\glie$ is finite dimensional and $v\in \tilde V^\prim$ is a primitive vector  that is killed by a large power of $\tilde f_k$ ($k=1,\dots ,r)$,  then this submodule is a finite dimensional  representation of $\glie$, which is irreducible in case $v$ is a weight vector of $\hlie$ (this weight is then necessarily dominant).
\end{lemma}
\begin{proof}
By the PBW-theorem,  the $\tilde\glie$-submodule of $\tilde V^\prim$ is also the
$\tilde\glie_-$-submodule generated by $\tilde V^\prim$.
In  $\tilde g$, the Serre element $\ad (\tilde e_k)^{1-c_{k,l}}\tilde e_l$ ($k\not= l$) commutes with every $\tilde f_k$ (see \cite{kac}, \S 3.3). Since it kills  $\tilde V^\prim$, it must be zero on the $\tilde\glie_-$-submodule generated by $\tilde V^\prim$. This proves the first assertion.

The second assertion  follows from Lemma 3.4 of \cite{kac}. The last assertion is then clear.
\end{proof}

In what follows $\lambda\in\hlie^*$ is a dominant  weight: for $k=1,\dots ,r$, $\lambda (\check{\alpha}_k)$ is a nonnegative real number. 
\\

\subsection*{Representations in $\hat\Bcal $} In what follows  we suppose our index set $\Ical$ endowed with a surjection $\pi :\Ical\to \{ 1,\dots ,r\}$ such that each fiber $\Ical_k:=\pi^{-1}(k)$ is countably infinite. We shall often write $\bar i$ for $\pi (i)$ and do likewise for the $\pi$-image of a sequence in $\Ical$. We will write $\Scal_\pt$ for $\Scal (\Ical_1)\times\cdots \times\Scal (\Ical_r)$

If $S$ is any sequence in $\{ 1,\dots ,r\}$, then we put
\[
\zeta (S):=\sum_{\bar I=S}\zeta_I ,
\]
where the sum is over all sequences in $\Ical$ that map under $\pi$ to $S$ (for  $S=\emptyset$, read $1$). The right hand side is an element of $\hat\Bcal$ that is invariant under the group 
$\Scal_\pt$. In fact, these elements give a basis of $\hat\Bcal ^{\Scal_\pt}$.

We take $c_{i,j}=c_{\bar i,\bar j}$ and $p_i:=\lambda (\check{\alpha}_{\bar i})$.
For this choice of coefficients, we put $\tilde f_k:=\sum_{i\in \Ical_k} \Phi_i$. We then have
\[
\tilde f_k \zeta (S)=\sum_{S=S''S'}(\lambda(\check{\alpha}_k)-c_{k,S'})\zeta (S''kS')
\]
Lemma \ref{lemma:efcommutator} suggests to put
\[
\tilde e_k \zeta(S):=
\begin{cases}
\zeta(S') & \text{ if $S=kS'$,}\\
0 & \text{otherwise}
\end{cases}
\]
(which makes $\tilde e_k$ independent of $\lambda$).
It is then clear that $\tilde e_k$ and $\tilde f_l$ commute when $k\not= l$ and that
$[\tilde e_k,\tilde f_k]$ multiplies $\zeta (S)$ by the scalar $\lambda (\check{\alpha}_k)-c_{k,S}$. Lemma \ref{lemma:residues} suggest that we have an  interpretation $\tilde e_k$ as a sum of  residues along divisors at infinity. There is a problem however since $\sum_{i\in \Ical_k} E_i$ does not make sense as a map defined on $\hat\Bcal$. Indeed, whereas 
$\tilde f_k$ makes sense on $\hat\Bcal$, there is no obvious way to define $\tilde e_k$ on that
space.  Yet Lemma \ref{lemma:efcommutator}  implies we may define it on $\hat\Bcal^{\Scal_\pt}$ as follows: 

\begin{lemma}\label{lemma:residuecrit}
Let $\zeta\in \hat\Bcal ^{\Scal_\pt}$.
Then for every $i\in \Ical_k$, we have
\[
(\tilde e_k(\zeta))_{\Ical-\{ i\}}=E_i(\zeta).
\]
In particular, $\tilde e_k(\zeta)=0$ if and only if $\zeta$ is regular along every hyperplane at infinity 
$(t_i=\infty)$ with $i\in \Ical_k$. 
\end{lemma}

If we combine this with \ref{lemma:flagresidue} and the antisymmetry of the Lie bracket, we obtain a way to express any Lie monomial in the  $\tilde e_k$'s as an iterated residue:

\begin{corollary}\label{cor:iteratedbracket}
Suppose we are in  the situation of  Lemma \ref{lemma:residuecrit}. Then for any sequence
$I=(i_N,\dots ,i_1)$  of length $N\ge 2$ in $\Ical$ without repetition we have for $\zeta\in \hat\Bcal^{\Scal_\pt}$
\[
\big([\cdots [(\tilde e_{\bar i_N},\tilde e_{\bar i_{N-1}}], \tilde e_{\bar i_{N-1}}],\cdots ,
 \tilde e_{\bar i_1}] \big)(\zeta)_{\Ical-\{ I\}}=
-\res_{(t_{\{ I\}}=\infty)}\res_I(\zeta).
\]
\end{corollary}

\begin{proposition}
The operators $\tilde e_k,\tilde f_k$, $k=1,\dots ,r$, define a representation of $\tilde\glie$ on $\hat\Bcal ^{\Scal_\pt}$ which satisfies for $k\not=l$ and $N\ge 1$ the identity  
\[
\frac{(\ad \tilde f_k)^N}{N!}\tilde f_l=(-c_{k,l}-1)\cdots (-c_{k,l}-(N-1))\Big(c_{k,l}\frac{\tau_k^N}{N!}\tilde f_l
+c_{l,k}\frac{\tau_l\tau_k^{N-1}}{(N-1)!}\tilde f_k\Big),
\]
so that indeed the Serre relations $\ad (\tilde f_k)^{1-c_{k,l}}\tilde f_l=0$ ($k\not=l$) are satisfied. 
\end{proposition}
\begin{proof}
Put $\check{\alpha}_k:=[\tilde e_k,\tilde f_k]$. We have seen that this operator is semisimple with integral eigenvalues ($\zeta(S)$ is an eigenvector with eigenvalue $\lambda(\check{\alpha}_k)-\sum_{i=1}^N c_{k,s_i}$). So  $\zeta(S)\in \hat\Bcal ^{\Scal_\pt}$ is an eigenvector of $\hlie$ of weight 
$\lambda-\sum_{i=1}^N \alpha_{s_i}$. 

The operator  $\tilde f_l$ changes the weight by $-\alpha_l$. Likewise, $\tilde e_l$ changes the weight by $\alpha_l$ and hence all the non-Serre relations are satisfied.
Corollary \ref{cor:serretype} shows  that  the displayed relations  also hold.
\end{proof}

When we regard  $\hat\Bcal ^{\Scal_\pt}$ as a $\tilde\glie$-module we shall denote it by 
$\tilde{V} (\lambda)$ and write $1_\lambda$ for its generator $1$. The $\hlie$-grading will be indicated by a subscript,  so that in the above proof, $\zeta(S)\in \tilde{V} (\lambda)_{\lambda-\sum_{i=1}^N \alpha_{s_i}}$. In particular, $\tilde{V} (\lambda)$ is a highest weight module of 
$\tilde \glie$ with highest weight $\lambda$. Notice that with $\check{\alpha}_k$ as above, we can now write 
\[
\tilde f_k \zeta(S)=
\sum_{S=S''S'}\check{\alpha}_k(\zeta(S))\zeta(S''kS').
\]
Denote by $V(\lambda)$ the $\tilde\glie_-$-submodule of $\tilde{V} (\lambda)$ generated by $1_\lambda$. It follows from the preceding proposition
that $V(\lambda)$ is then also invariant under $\tilde\glie_+$, hence is a $\tilde \glie$-module. \\

\begin{remark}
Lemma \ref{lemma:residuecrit} characterizes the primitive part $\tilde{V} (\lambda)^\prim$ of
$\tilde{V} (\lambda)$ as the subspace of polydifferentials that are regular at the hyperplanes $t_i=\infty$. This is independent of $\lambda$ (and we may even define $\hat\Bcal^\prim$ as the space of such polydifferentials, although  we did not define an action of $\tilde\glie$ on $\hat\Bcal$).

Clearly the generator  $1_\lambda\in\tilde{V} (\lambda)$ is primitive.
\end{remark}

\begin{theorem}\label{thm:grep}
The Lie algebra $\tilde\glie$ leaves invariant the $\tilde\glie_-$-submodule of $\tilde{V} (\lambda)$ generated  by the primitive subspace $\tilde{V} (\lambda)^{\prim}$ and  acts on that space through $\glie$. In particular, the subrepresentation $V(\lambda)$ 
generated by $1_\lambda$ is a highest weight representation of  $\glie$. If $\lambda$ is an integral weight, then this representation is integrable in the sense that each of the $e_k$ and $f_k$ acts on it in a  locally nilpotent fashion. In case the given Cartan matrix is that of a finite dimensional Lie algebra, then 
$V(\lambda)$ is finite dimensional and irreducible.
\end{theorem}
\begin{proof} This is a direct application of Lemma \ref{lemma:primg}, where for the last half we invoke  Corollary \ref{cor:power}.
\end{proof}

Notice that we do not claim that $\tilde{V} (\lambda)$ is a representation of $\glie$. Indeed, it is not true in general that $\ad (\tilde e_k)^{1-c_{k,l}}\tilde e_l$ vanishes on that space. 

\begin{example}
Assume  that $r=1$. Then $\tilde\glie=\glie\cong\slin (2)$  and by Corollary \ref{cor:polynomial}, $\hat\Bcal^{\Scal}$ is the polynomial algebra on the generator $\zeta$ defined by $\zeta (z)=\sum_{i\in\Ical}(t_i-z)^{-1}dt_i$. One verifies that the operator $e$ is simply derivation: $e(\zeta^N)=N\zeta^{N-1}$. Let us identify the weight $\lambda$ that turns $\hat\Bcal^{\Scal}=\CC[\zeta]$ into the   $\slin (2)$ representation $\tilde V(\lambda)$ with its value on the unique simple coroot. Then another straightforward  computation shows that $f$ then sends $\zeta^N$ to $(\lambda -2N)\zeta^{N+1}$.
 \end{example}

The polar divisor of any $\zeta\in \hat\Bcal^{\Scal_\pt}$ is in general much smaller than that of an arbitrary member of 
$\hat\Bcal$. For instance, if  $i,j\in\Ical$ are distinct, but such that $\bar i=\bar j$, then $\zeta$ has no poles along the diagonal hyperplane $t_i=t_j$. To see this observe that $\tilde\zeta:=(t_i-t_j)\zeta$ has no pole along  $t_i=t_j$ and that since $\zeta$ is $\Scal_\pt$-invariant, interchanging $t_i$ and $t_j$ turns $\tilde\zeta$ into $-\tilde\zeta$. So $\tilde\zeta$ is divisible by $(t_i-t_j)$ and hence $\zeta$ has no pole along $t_i=t_j$.

We can do even better  on $V(\lambda)$:

\begin{lemma}\label{lemma:diagonal2} 
Let $I$ be a finite  sequence in $\Ical$. If   $\res_{I}$ is nonzero on $V(\lambda)$, then 
for every initial part $J$ of $I$, $\alpha_{\{\bar J\}}$ is a root. 
\end{lemma}
\begin{proof} Let $I=(i_N,\dots ,i_1)$ and let $\zeta\in V(\lambda)$ be homogeneous, of degree $m$, say.  If $N>m$, then 
$\res_I(\zeta)  =0$ and so there is nothing to prove. We proceed with downward induction on $N$.

If for some $i\in \Ical$,  $\res_{Ii}\zeta\not=0$, then we may apply our induction hypothesis and conclude that for every initial part $J$ of $Ii$, $\alpha_{\{\bar J\}}$ is a root. It remains to deal with the case when $\res_{Ii}\zeta$ vanishes for all $i$. Since $\res_{Ii}\zeta=\res_{(t_{i_1}\to t_i)}\res_{I}\zeta$, this means that 
 the poles of $\res_{I}(\zeta)$ that involve the coordinate $t_{\{ I\}}$ can only occur for $t_{\{ I\}}=z$ and $t_{\{ I\}}=\infty$ and so $-\res_{(t_{\{ I\}}=\infty)}\res_{I}\zeta=\res_{(t_{\{ I\}}=z)}\res_{I}\zeta$ by the residue theorem. Recall that  $\res_{(t_{\{ I\}}=\infty)}\res_{I}\zeta =(-e_{\bar I}\zeta)_{\Ical-\{ I\}}$ and so if this is nonzero, then  for every initial part $J$ of $I$, $e_{\bar J}\not=0$ and hence 
 $\alpha_{\{\bar J\}}$ is a root. Otherwise, $\res_{I}\zeta$ is regular for a generic value of the coordinates $t_j$, $j\in \Ical$. But a Riemann sphere has no nonzero holomorphic differential and since $\res_{I}\zeta$ involves  $dt_{\{ I\}}$, this must imply that $\res_{I}\zeta$ is identically zero.
\end{proof}

We recall that if $\glie$ is simple and  finite dimensional, then there is a unique \emph{highest root} $\tilde\alpha$ relative to the root basis $(\alpha_1,\dots ,\alpha_r)$. It is the unique long root that also dominant and also characterized by the property that  for no $k=1,\dots, r$,
$\tilde\alpha +\alpha_k$ is a root.

\begin{corollary}\label{cor:diagonal2}
Suppose  that $\glie$ is simple and  finite dimensional.
If $I$ is a finite sequence in $\Ical$ such that $\alpha_{\{\bar I\}}$ is the highest root, then for every $\zeta\in V(\lambda)$,  the poles of $\res_{I}(\zeta)$ that involve the coordinate $t_{\{ I\}}$ can only occur for  $t_{\{ I\}}=z$ and $t_{\{ I\}}=\infty$ and we have  $(e_{\tilde\alpha}\zeta)_{\Ical-\{I\}}=\res_{(t_{\{ I\}}=z)}\res_I (\zeta)$. 
\end{corollary}

The automorphism group $\aut (\Pbold)$ of $\Pbold$ acts on $\Pbold^{\Ical}$. In fact, it acts on the projection $\Pbold^{\Ical}\times\Pbold\to\Pbold$. This action preserves the divisors of the form $(t_i=t_j)$, and $(t_i=z)$, but not the divisors $(t_i=\infty)$. The stabilizer of $\infty$ in $\aut (\Pbold)$ does have that property however and indeed, it leaves every $\zeta_I$ invariant and hence acts as the identity on $\Bcal$. We should not expect this action to happen for all of  $\aut (\Pbold)$. Indeed, if $\sigma=\binom{a\, b}{c\, d}\in \SL (2,\CC)$, then a straightforward computation shows that
\[
\sigma^*\zeta_I (z)=\frac{cz+d}{ct_{b(I)}+d}\zeta_I(z).
\]
According to Lemma \ref{lemma:mixedshuffle} the left hand side equals (at least for $c\not=0$)
\[
\sum_{I=I''I'} (-1)^{|I''|}\zeta_{I''{}^*}(-d/c).\zeta_{I'}(z).
\]
So on $\hat\Bcal^{\Scal_\pt}$ we find
\begin{multline*}
\sigma^*\zeta (S)=\sum_{S=S''S'} (-1)^{|S''|}\zeta (S''{}^*)(-d/c). \zeta (S')=\\
=\sum_{T} (-1)^{|T|}\zeta (T)(-d/c). \tilde e_{T}(\zeta (S)),
\end{multline*}
where the sum is over all finite sequences $T$ in $\{ 1,\dots ,r\}$ (but we get only a nonzero contribution from $T$ if it appears as the initial part of $S$). To sum up, on
$\hat\Bcal^{\Scal_\pt}$ we have that
\[
\sigma^*=\sum_{T} (-1)^{|T|}\zeta (T)(-d/c). \tilde e_{T}.
\]
On  the primitive part of $\tilde V (\lambda)^\prim$ this reduces to the term corresponding to $T=\emptyset$, which is just $\zeta (S)$.  So we find:

\begin{corollary}\label{cor:slinvariance}
The primitive part of $\hat\Bcal^{\Scal_\pt}$ is left pointwise fixed under the action of  
the automorphism group of $\Pbold$. 
\end{corollary}

In the next section we shall also need to know the infinitesimal (right) action of $\binom{0\, 0}{1\, 0}\in \slin (2,\CC)$. 
 on $\hat\Bcal^{\Scal_\pt}$. A similar argument shows that
\[
\begin{pmatrix}
0& 0\\1& 0
\end{pmatrix}^*=\sum_{T} (-1)^{|T|}\omega (T). \tilde e_{T},
\] 
where
$\omega (T)=\sum_{\bar I=T} \omega_I$.

\subsection*{Tensor products}
We generalize the above to the case of a tensor product of highest weight representations.
We fix an $n$-tuple of dominant weights, $\lambdabold=(\lambda_1,\dots ,\lambda_n)$ and
instead of working with the base $\Abold$, we use the base $\Abold^n$: rather than a single variable $z$ we have $n$ variables $\zz=(z_1,\dots ,z_n)$. 
For an $n$-tuple $\Ibold=(I^1,\dots ,I^n)$ of  sequences in $\Ical$ we consider the relative polydifferential 
\[
\zeta_{\Ibold }: \zz=(z_1,\dots ,z_n)\mapsto\zeta_{\Ibold}(\zz)=\zeta_{I^1}(z_1)\cdot\zeta_{I^2}(z_2)\cdot\dots \cdot\zeta_{I^n}(z_n).
\]
It clear that this polydifferential vanishes unless  the sequence $I^1\cdots I^n$ obtained by juxtaposition is without repetition. We denote by $\Bcal_{n}$ be the graded  vector space spanned by these polydifferentials, by  $\hat\Bcal_n^d$ the completion of $\Bcal_n^d$ of (form) degree $d$ which allows for infinite sums of these polydifferentials and put 
$\hat\Bcal_n:=\oplus_d\hat\Bcal_n^d$.
It may be worthwhile to observe that $\zeta_{\Ibold }$ is invariant under the stabilizer of $\infty$ in the automorphism group $\Pbold$, in other words, under the automorphism group of $\Abold$.

Given an $n$-tuple $\Sbold=(S^1,\dots ,S^n)$ of  sequences in $\{1,\dots ,r\}$, we observe that
\[
\zeta (\Sbold)(\zbold):=\sum_{\overline{\Ibold}=\Sbold} \zeta_{\Ibold}(\zbold)=\zeta (S_1)(z_1)\cdots 
\zeta (S_n)(z_n),
\]
where the sum is over all $n$-tuples of sequences $\Ibold=(I^1,\dots ,I^n)$ in $\Ical$ whose juxtaposition is without repetition and map under $\pi$ to $\Sbold$. These elements form a $\CC$-basis of $\hat\Bcal_n^{\Scal_\pt}$ and so the above factorization defines an isomorphism 
\[
\hat\Bcal _{n}^{\Scal_\pt}\cong \hat\Bcal^{\Scal_\pt}\otimes_\CC\cdots \otimes_\CC\hat\Bcal^{\Scal_\pt}.
\]

\begin{remark}\label{ex:sltwo}
Assume that $r=1$, so that  $\glie=\slin (2)$. It then follows from Corollary \ref{cor:polynomial} that $\hat\Bcal _{n}^{\Scal}$ is a polynomial algebra on the $n$ generators $\sum _{i\in\Ical} (t_i-z_\nu)^{-1}dt_i$. This yields  the free $\CC$-basis
\[
\zeta (1^{k_1},\dots ,1^{k_n})=\sum_{\{(X_\pt), | X_\pt|=k_\pt\}} \prod_{\nu =1}^n\prod_{x\in X_\nu}\frac{dt_x}{t_x-z_\nu},
\]
where the sum is over $n$-tuples of disjoint subsets $(X_1,\dots ,X_n)$ of $\Ical$ with $| X_\nu|=k_\nu$.
\end{remark}

The action of $\tilde f_k$ operating on the $\nu$th factor with dominant weight $\lambda ^{(\nu)}$ is denoted  
$\tilde f_k^{(\nu)}$. The sum $\sum_{\nu =1}^n \tilde f_k^{(\nu)}$ acts as  $\tilde f_k$ in the tensor representation and hence is simply denoted $\tilde f_k$. Notice  however that we can define
$\tilde f_k$ without reference to the tensor decomposition above as
\[
\tilde f_k=\sum_{i\in\Ical_k}dt_i\Big(\sum_{\nu=1}^n \frac{\lambda_\nu(\check{\alpha}_k)}{t_i-z_\nu}-
\sum_{j\in\Ical, j\not=i} \frac{c_{k, \bar j}}{t_i-t_j}\iota_{\partial/\partial t_j}\Big).
\]
For the tensor action of $\tilde e_k$ (as $\sum_{\nu =1}^n \tilde e_k^{(\nu)}$) the situation is even better, for Lemma \ref{lemma:residuecrit} remains valid in this multivariable setting:

\begin{lemma}\label{lemma:residuecrit2} 
If $i\in\Ical_k$ and $\zeta\in \hat\Bcal _{n}^{\Scal_\pt}$, then
\[
\tilde e_k (\zeta)_{\Ical-\{ i\}}=E_i(\zeta).
\]
\end{lemma}
\begin{proof}
It is enough to verify this in case $\zeta=\zeta (\Sbold)$, where $\Sbold =(S^1,\dots ,S^n)$ is an $n$-tuple of sequences in $\{1,\dots ,r\}$. If 
$\Ibold=(I^1,\dots ,I^n)$ is an $n$-tuple of sequences in $\Ical$ whose concatenation is without repetition and which maps under $\pi$ to $\Sbold$, then the value of $E_i$ on
$\zeta_{\Ibold}$ is zero unless $i$ is the first term of some $I^\nu$: $I^\nu=iI'{}^\nu$, in which case we get $\zeta_{(I^1,\dots , I'{}^\nu,\dots ,I^n)}$. So if we take the sum over such 
$\Ibold$ we find the $(\Ical-\{ i\})$-component of $\tilde e_k(\zeta (\Sbold))$.
\end{proof}

This completely describes $\hat\Bcal _{n}^{\Scal_\pt}$ as a tensor product  of representations of $\tilde\glie$. We will denote it by $\tilde{V} (\lambdabold)$ and write $1_{\lambdabold}$ for its generator $1$. 
We obtain the following generalization of Theorem \ref{thm:grep}.

\begin{theorem}\label{thm:grep2}
We have a natural identification 
\[
\tilde{V} (\lambdabold)\cong\tilde{V} (\lambda_1)\otimes_\CC\cdots \otimes_\CC
\tilde{V} (\lambda_n).
\]
The primitive subspace $\tilde{V} (\lambdabold)^\prim\subset\tilde{V} (\lambdabold)$ is  the subspace  consisting of forms that are regular at every hyperplane at infinity $(t_i=\infty)$, $i\in \Ical$, and the  $\tilde\glie_-$-submodule of  $\tilde{V} (\lambdabold)$ it generates   is acted on by the Lie algebra $\tilde\glie$  via $\glie$ (so it is in fact a  $\tilde\glie_-$-module). In particular, 
\[
V(\lambdabold):= V(\lambda_1)\otimes_\CC\cdots \otimes_\CC V(\lambda_n)
\]
is the smallest subspace of $\hat\Bcal _{n}^{\Scal_\pt}$ that contains $1_{\lambdabold}$ and is invariant under the operators  $f_k^{(\nu)}$ and $e_k^{(\nu)}$; it is the tensor product of $n$ highest weight representations of $\glie$. It is integrable if all the $\lambda_k$'s are integral.
\end{theorem}

We next state two important properties of the elements of $V(\lambdabold)$ and 
$V(\lambdabold)^\prim$. First, Lemma  \ref{lemma:diagonal2} almost immediately generalizes to:

\begin{theorem}\label{thm:diagonaltensor}
If $I=(i_N,\dots ,i_1)$ is a  sequence in $\Ical$ such that  $\res_{I}$ is nonzero on $V(\lambdabold)$, then  for every initial part $J$ of $I$, $\alpha_{\{\bar J\}}$ is a root. 

In case the last element $\alpha_{\{ \bar I\}}$ is the highest root $\tilde\alpha$, then for  every $\zeta\in V(\lambda)$,  the poles of $\res_{I}(\zeta)$ that involve the coordinate $t_{\{ I\}}$ can only occur where $t_{\{ I\}}$ takes a value in $\{\infty, z_1,\dots z_n\}$, where we may omit $\infty$ in case  $\zeta\in V(\lambdabold)^\prim$. We also have   $(e_{\tilde\alpha}\zeta)_{\Ical-\{I\}}=\sum_{\nu=1}^n\res_{(t_{\{ I\}}= z_\nu)}\res_I (\zeta)$.
\end{theorem}

We also have the obvious extension of Corollary \ref{cor:slinvariance}. 

\begin{theorem}\label{thm:slinvariances}
The primitive part of $\hat\Bcal_n^{\Scal_\pt}$ is left pointwise fixed under the action of  
the automorphism group of $\Pbold$. 
\end{theorem}
\begin{proof}
We only need to verify this infinitesimally: that any primitive element of $\hat\Bcal_n^{\Scal_\pt}$ is killed
by the Lie algebra of $\aut (\Pbold)$. This is clearly so for all the standard generators of this Lie algebra except $\binom{0\, 0}{1\, 0}\in \slin (2,\CC)$. Its action is however on a given tensor factor given by the expression 
$\sum_{T} (-1)^{|T|}\omega (T). \tilde e_{T}$. Hence the same is true on the full tensor product. The last expression clearly vanishes on  $\hat\Bcal^{\Scal_\pt}$.
\end{proof}

Note: \emph{From now on we assume the weights $\lambda_1,\dots ,\lambda_n$ to be integral.}
\\ 

For the KZ-equation we shall have  to consider the subspace $V(\lambdabold)^\glie$ of $\glie$-invariants in $V(\lambdabold)$. This is just  $V(\lambdabold)_0^\prim$, the primitive part of $V(\lambdabold)$ of weight zero.  (Note that $V(\lambdabold)_0\not=0$ implies that $\sum_\nu \lambda_\nu$ is a sum of positive roots.) 

We invoke the representation theory of $\slin (2)$ to deduce:

\begin{proposition}
The intersection of $\sum_k \tilde f_k \tilde V(\lambdabold)$ with $V(\lambdabold)_0$ is
the subspace $(\sum_k \tilde f_k \tilde V(\lambdabold))_0$. In particular,  $V(\lambdabold)^\glie$ embeds in $\tilde{V} (\lambdabold)_{0,\coprim}$.
\end{proposition}
\begin{proof}
The grading of $\tilde V(\lambdabold)$ by the weights of $\hlie$ shows that we have 
$\sum_k \tilde f_k \tilde V(\lambdabold)\cap V(\lambdabold)_0=\sum_k \tilde f_k \tilde V(\lambdabold)_{\alpha_k}$. So it suffices to show that  $\tilde f_k \tilde V(\lambdabold)_{\alpha_k}\cap 
V(\lambdabold)=\tilde f_k V(\lambdabold)_{\alpha_k}$ for every $k$. This makes it an issue about  
$\slin (2)$. With induction it is easily shown that for every $v\in \tilde V(\lambdabold)_{\alpha_k}$ with
$\tilde f_k v\in V(\lambdabold)$ we have  
\[
v\equiv \frac{(-1)^p}{p! (p+1)!}\tilde f_k^p\tilde e_k^p (v)\pmod{V(\lambdabold)}.
\]
Since $\tilde e_k^p (v)=0$ for $p$ large, the claim follows.
\end{proof}

\subsection*{The passage  to differential forms}\label{subsect:diffforms}
In what follows we assume $n\ge 2$. We take as our base variety $U_n$  the set of  $(z_1,\dots ,z_n)\in\Abold^n$ with $z_1,\dots ,z_n$ pairwise distinct and  summing up to  $0\in\Abold$; in other words $U_n$ is the  standard arrangement complement of type $A_{n-1}$. It is better however to refrain from choosing an origin for $\Abold$ and to think of $U_n$ in modular terms:  if we endow $\Abold$ with the constant differential $dz$, then we easily see that  $U_n$ may be identified with the moduli space of smooth genus zero curves $C$ endowed with an embedding of  $\{\infty ;1,\dots ,n\}$  in $C$ and with a differential $dz$ on the complement of the image of $\infty$ that is invariant under the automorphism group of that complement. The latter is also equivalent to the choice of a nonzero tangent vector of $C$ at the image of $\infty$, so this makes $U_n$ a (trivial)  $\CC^\times$-bundle over $\Mcal_{0,1+n}$.

We now assume that $V(\lambdabold)_0\not=0$, so that  $\sum_\nu\lambda_\nu$ is a sum of simple roots: $\sum_\nu\lambda_\nu=\sum_{k=1}^r m_k\alpha_k$ with $m_k\in \ZZ_{\ge 0}$. We put $\mbold:=(m_1,\dots ,m_k)$ and $m:=\sum_k m_k$, so that  $V(\lambdabold)_0$ lies in the homogeneous summand of multi-degree $\mbold$.
Let $M\subset \Ical$ be a finite subset so that  $M_k:=M\cap \Ical_k$ has cardinality $m_k$. 
Our use of the symbol $M$ implies that this decomposition is understood.
We denote by $\orient_M:=\wedge^m \ZZ^M$ the sign representation of $\Scal (M)$ and we
put $\Scal (M_\pt):=\Scal_\pt\cap\Scal (M)= \Scal(M_1)\times\cdots\times \Scal (M_r)$.

We denote  by  $\Bcal_{n,M}$ the corresponding graded  algebra of relative  polydifferentials on $\Pbold^{M}_{U_n}
=(\Pbold^{M_1}\times\cdots \times\Pbold^{M_r}\times U_n)/U_n $. We regard this as a subalgebra of  $\Bcal _{n}$ via pull-back.  It is multigraded by $r$-tuples of nonnegative integers and has $\mbold$ as highest  multi-degree. 
Consider the homogeneous elements of  $\Bcal_{n,M}^{\Scal(M_\pt)}$ defined by 
\[
\xi_k:=\sum_{i\in M_k} \Big(\sum_{\nu=1}^n \frac{\lambda_\nu(\check{\alpha}_k)}{t_i-z_\nu}
-\sum_{j\in M -\{ i\}} \frac{c_{k,\bar j}}{t_i-t_j}\Big) dt_i.
\]

\begin{lemma}\label{lemma:finitepart}
The map which assigns to an element of $\tilde{V} (\lambdabold)$ the sum of its $X$-components, where $X$ runs over the subsets of $ M$,
maps onto $\Bcal_{n,M}^{\Scal( M_\pt)}$ and so identifies the latter with $\tilde{V}_M(\lambdabold)$. It is an isomorphism in nonnegative weights: $\tilde{V} (\lambdabold)_{\ge 0}\cong \tilde{V}_M(\lambdabold)$.  The transferred action of $\tilde e_k$ to $\tilde{V}_M(\lambdabold)$ is the obvious one 
(and given by residues as in Lemma \ref{lemma:residuecrit2}) and the same is true for  $\tilde f_k$  on the summands of multi-degree strictly lower than $\mbold$. 

Moreover, if $a_1,\dots ,a_k$ are nonzero complex numbers and $m_k\ge 1$ for all $k$, then  the $M$-component of $\tilde{V} (\lambdabold)_{0,\coprim}$,  $\tilde{V} (\lambdabold)^M_{0,\coprim}$, gets identified with that of
$\tilde{V}_M(\lambdabold)/((\sum_ka_k\xi_k)\tilde{V}_M(\lambdabold))$.
So this embeds ${V}(\lambdabold)^\glie$ in $\tilde{V} (\lambdabold)^M_{0,\coprim}$.
\end{lemma}
\begin{proof}
All but the last of these assertions follow from the fact that every $\alpha\in \Bcal_{n,M}^{\Scal_\pt}$ of degree $d$ is uniquely written as $\sum_{X\in\Ical, |X|=d} \pi_X^*\alpha_X$, where $\alpha_X$ is a rational polydifferential on $\Pbold^X_{\Abold^n}$.
To prove the last one, let us first  observe that if $\zeta=\prod_{j\in M} dt_j$, then 
for $k=1,\dots ,r$, then every element of $\Bcal_{n,M}$ of degree  $m-1$ is a linear combination
of the forms $\iota_{\p/\p t_i}\zeta$ with $i\in M$ with rational functions as coefficients.
Now if $k\in\{ 1,\dots ,r\}$ and $i\in M_k$, then 
\[
(\sum_{k=l}^ra_l\xi_l)\iota_{\p/\p t_i}\zeta=a_k\Big(\sum_{\nu=1}^n \frac{\lambda_\nu (\check{\alpha}_k)}{t_i-z_\nu}
-\sum_{j\in M -\{ i\}} \frac{c_{\bar i,\bar j}}{t_i-t_j}\Big)\zeta
\]
is also the $ M$-component of $\tilde f_k(\iota_{\p/\p t_i}\zeta)$. The last statement now follows easily.
\end{proof}

The moduli space of injections of the disjoint union of $\{1,\dots ,n\}\sqcup M$ in $\Abold$ given up to translations is also the moduli space of triples 
\[
(C, z\sqcup t: \{\infty ,1,\dots ,n\}\sqcup M\hookrightarrow C, dz),
\] 
where $C$ is a complete smooth curve of genus zero, $z\sqcup t$ is an embedding and 
$dz$ is a nonzero differential on $C$ invariant under $\aut (C, z(\infty))$. We denote it by 
$U_{n,M}$. Ignoring the embedding of $M$  defines an evident morphism $U_{n,M}\to U_n$.  We make  this morphism proper by means of a relative Deligne-Mumford-Knudsen compactification $U^+_{n,M}\to  U_n$. Here $U^+_{n,M}$ is a $\CC^\times$-bundle over a moduli space of  stable  pointed genus zero curves: it parameterizes triples   $(C, z\sqcup t, dz)$ as before, where we now allow  the pair $(C,z\sqcup t)$ to be a stable pointed curve, but insist that if we ignore the embedding of $M$ and contract the irreducible components of $C$ as to make it stable, the result is an $(1+n)$-pointed curve that is   \emph{smooth}. In particular, we get a retraction of $C$ onto a distinguished component $C_o$ such that its composite with $z$ is injective.  This component $C_o$ must then contain all but at most one of the images of $z$ and it is on this component that we assume $dz$ is defined.
This exhibits the desired morphism $p^+:U^+_{n,M}\to U_n$. It is proper, indeed. In fact,
if we extend the definition of the reduced effective divisors $D_X$ in an obvious manner as (relative) divisors on $\Pbold^M_{U_n}$, or rather $\Pbold^M_{\Pbold^n}$:
\begin{align*}
D^f_{n,M}:=&\sum_{i\in M}\sum_{\nu=1}^n (t_i=z_\nu) +\sum_{\{i\not=j\}\subset M}(t_i=t_j),\\
D_{n,M}:=& D^f_{n,M} +\sum_{i\in M}(t_i=\infty).
\end{align*}
(the superscript in $D^f_{n,M}$ stands for finite), then $U^+_{n,M}$ is obtained from $\Pbold^M\times U_n$ by a blowing up process that is minimal for the property of turning
$D_{n,M}$ into a normal crossing divisor $\Delta_{n,M}= U^+_{n,M}-U_{n,M}$. The generic point of an irreducible component of this divisor parameterizes one point unions of two smooth rational curves with the disjoint union of $\{ \infty,1,\dots ,n\}$ and $M$ embedded in its smooth part such that besides obeying the usual stability condition (every connected component of the smooth part contains at least two of these points) we have that one of the two connected components meets $\{ \infty,1,\dots ,n\}$ either in a singleton or not at all. We denote this irreducible component of the boundary divisor accordingly as $\Delta_\infty(X) $, $\Delta_\nu (X)$ or $\Delta (X)$,  where $\nu=1,\dots ,n$ and  $X\subset M$ is nonempty and has at least two distinct elements in the last case. In terms of  the configuration space of maps $X\to \Pbold$ this corresponds to a confluence of the members $X$ (where in the first case resp.\ second case the confluence is towards $\infty$ resp.\ $z_\nu$).  We denote by $\Delta^f_{n,M}$ the `finite' part of $\Delta_{n,M}$, that is, the sum of the divisors $\Delta_\nu (X)$  and  $\Delta (X)$. This is indeed the full preimage of $D^f_{n,M}$ under the blowup. It follows from Corollary \ref{cor:intrinsic} that we have a natural identification $\CC[U_n]\otimes\Bcal_{n,M}\cong \bigoplus_{X\subset  M} H^0(\Pbold^X_{U_n}, \Omega_{\Pbold^X_{U_n}}^{|X|}(\log D_{n,X}))\otimes\orient(X)$.
This gives rise to an isomorphism of $\Ocal_{U_n}$-modules:
\[
\Ocal_{U_n}\otimes\Bcal^m_{n,M}\cong
p^+_*\Omega^m_{U^+_{n,M}/U_n}(\log \Delta_{n,M}))\otimes\orient(M).
\]
If we combine this with Corollary \ref{cor:intrinsic} and  the Lemmas
\ref{lemma:residuecrit2}   and \ref{lemma:finitepart}, we find:

\begin{proposition}\label{prop:findim}
We have a natural identification of  $\tilde{V} (\lambdabold)_{0}$ with the subspace of 
$H^0(U_{n,M}, 
\Omega_{U_{n,M}/U_n}^{m}(\log \Delta_{ M,n}))\otimes_{\Scal (M_\pt)} \orient (M)$
of relative logarithmic $m$-forms that vanish along the hyperplanes $z_\nu=\infty$, $\nu=1,\dots ,n$.
This restricts to
an isomorphism of $\tilde{V} (\lambdabold)_0^\prim$ with the  corresponding subspace 
of $H^0(U_{n,M}, 
\Omega_{U_{n,M}/U_n}^{m}(\log \Delta^f_{ M,n}))\otimes_{\Scal (M_\pt)} \orient (M)$ and yields trivializations of bundles over $U_n$:
\begin{align*}
\Ocal_{U_n}\otimes_\CC \tilde{V} (\lambdabold)_0
&\cong (p^+_*\Omega_{U^+_{n,M}/U_n}^{m}(\log \Delta_{ M,n}))\otimes_{\Scal (M_\pt)} \orient (M),\\
\Ocal_{U_n}\otimes_\CC\tilde{V} (\lambdabold)^{\prim}_0
&\cong (p^+_*\Omega_{U^+_{n,M}/U_n}^{m}(\log \Delta^f_{ M,n}))\otimes_{\Scal (M_\pt)} \orient (M).\\
\end{align*}
In particular, $\Ocal_{U_n}\otimes_\CC\tilde{V} (\lambdabold)^{\glie}$ embeds in the last module.
\end{proposition}

\section{Identification of the KZ connection} 

We continue with the situation of the previous section. So we
have the $n$-tuple of integral dominant weights $\lambdabold=(\lambda_1,\dots ,\lambda_n)$
and regard $V(\lambdabold)=V(\lambda_1)\otimes\cdots\otimes 
\tilde{V}(\lambda_n)$ as a representation of $\glie$. We let $\mbold =(m_1,\dots ,m_r)$ and  the finite subset $M\subset\Ical$ be as in Section \ref{sect:reps}.

\subsection*{The KZ-connection}\label{ssect:differentpres}
The Knizhnik-Zamolodchikov connection requires the choice of a  \emph{Casimir element}, that is, a symmetric tensor  $C\in \glie\otimes\glie$ that is invariant under $\glie$ (acting adjointly on each factor). 

We  denote by  $q^C$ the quadratic form on $\glie^*$ attached to $C$: $q^C(a)=\half C(a,a)$ (note the factor $\half$). We do not need to assume that $C$ is nondegenerate, but we do suppose that $q^C(\alpha_k)\not=0$ for every simple root $\alpha_k$. This implies that the generalized Cartan matrix is symmetrizable and that $C$ is nondegenerate on $(\glie/\hlie)^*\subset\glie^*$. 

We use this tensor $C$ to obtain a slightly different presentation of $\tilde\glie$: if $\alpha_k$ is the positive root attached to $e_k$ (so $[\check{\alpha}_l,f_k]=-\alpha_k(\check{\alpha}_l)f_k$), then we replace $f_k$ by $\hat f_k:=q^C(\alpha_k)f_k=\half C(\alpha_k,\alpha_k)f_k$.
We retain $e_k$ and so $\hat h_k:=[e_k,\hat f_k]=\half C(\alpha_k,\alpha_k)\check{\alpha}_k$ has now  the property that $\lambda (\hat h_k)= C(\lambda ,\alpha_k)$ for every $\lambda\in\hlie^*$. Then $\tilde\glie$ is presented in terms of the symmetric matrix $(C(\alpha_k,\alpha_l))_{k,l}$: 
\[
[\hat h_k,e_l]=C(\alpha_k,\alpha_l)e_l,\quad  [\hat h_k,\hat f_l]=-C(\alpha_k,\alpha_l)
\hat f_l,\quad   [\hat h_k,\hat h_l]=0.
\]
We still need to impose the Serre relations (which involve the possibly nonsymmetric Cartan matrix), but these are just the ones that make $C$ nondegenerate on the subspace
$(\glie/\hlie)^*\subset\glie^*$.   In this setup  $\lambda(\check{\alpha}_k)$ becomes  $C(\lambda,\alpha_k)$, $c_{k,l}$ becomes $C(\alpha_k,\alpha_l)$, and $\Phi_k$ is replaced by
\[
\hat\Phi_k :=\sum_{i\in\Ical_k}dt_i\Big(\sum_\nu \frac{C(\lambda_\nu,\alpha_k)}{t_i-z_\nu}-\sum_{j\not=i} C(\alpha_k,\alpha_{\bar j})\frac{dt_j}{t_i-t_j}\iota_{\p/\p t_j}\Big).
\]
For $1\le \nu<\mu\le n$, let $C_{\nu,\mu}$ be the endomorphism of $V(\lambdabold)$ obtained by letting  $C$ act trough the tensor factors indexed by  $\nu$ and $\mu$. This operator commutes with with the diagonal action of $\glie$ and hence preserves the $\glie$-isotypical summands.
Then the  corresponding KZ connection $\nabla^C_{KZ}$ on $\Ocal_{U_n}\otimes_\CC V(\lambdabold)$ is defined by the $\End (V(\lambdabold))$-valued differential
\[
A^C_{KZ}:=\sum_{1\le \nu<\mu\le n} C_{\nu,\mu}\frac{d(z_\nu-z_\mu)}{z_\nu-z_\mu}.
\]
This is a connection  with logarithmic singularities at infinity. It is easily shown to be  flat so that we get a local system 
$\KZ^C(\lambdabold)\subset \Ocal_{U_n}\otimes_\CC V(\lambdabold)$.

\begin{remark}[Comparison with the SV-map]
This is essentially the situation considered by  Schechtman-Varchenko in \cite{sv} from the outset.
Our space $\tilde{V}_M(\lambdabold)$ is basically the one they construct for the case of a symmetrizable 
Cartan matrix. They identify the action of the operators $\hat f_1,\dots ,\hat f_r$, but there are no operators 
$\tilde e_1,\dots, \tilde e_r$ acting. So the coprimitive quotients can (and do) appear there, but primitive subspaces cannot.    
\end{remark}

\subsection*{A local system of rank one}
Central in the subsequent discussion will be the following differential associated to $C$ (a
formal expression, for this is an infinite sum):
\begin{multline*}
\eta _{\lambdabold}^C:=\sum_{\nu=1}^n\sum_{i\in \Ical} C(\alpha_{\bar  i},\lambda_\nu)\frac{d(t_i-z_\nu)}{t_i-z_\nu}\\
-\half\sum_{i,j\in\Ical, i\not=j} C (\alpha_{\bar i},\alpha_{\bar j})\frac{d(t_i-t_j)}{t_i-t_j}
-\sum_{1\le \nu<\mu\le n}C(\lambda_\nu,\lambda_\mu)\frac{d(z_\nu-z_\mu)}{z_\nu-z_\mu}.
\end{multline*}
The corresponding relative form $\xi_{\lambdabold}^C:=(\eta_{\lambdabold}^C)_{\text{rel}}$ is obtained by ignoring the
$dz_\nu$-terms. Consider the finite subsums $\eta_{\lambdabold,M}^C$ and $\xi_{\lambdabold,M}^C$ that involve the factors indexed by $M$. So 
\begin{multline*}
\xi^C_{\lambdabold,M}:=(\eta_{\lambdabold,M}^C)_{\text{rel}}=
\sum_{i\in M} \Big(\sum_{\nu=1}^n \frac{C(\alpha_{\bar  i},\lambda_\nu)}{t_i-z_\nu}
-\sum_{j\in M -\{ i\}} \frac{C (\alpha_{\bar i},\alpha_{\bar j})}{t_i-t_j}\Big) dt_i\\
=\sum_{k=1}^r q^C(\alpha_k)\sum_{i\in M_k} \Big(
\sum_{\nu=1}^n \frac{\lambda_\nu(\check{\alpha}_k)}{t_i-z_\nu}-\sum_{j\in M -\{ i\}} 
\frac{ c_{k,\bar j} }{t_i-t_j}\Big) dt_i
=\sum_{k=1}^r q^C(\alpha_k)\xi_k.
\end{multline*}
Since  each $q^C(\alpha_k)$ is nonzero, this is an element of 
$\tilde{V} (\lambdabold)^{ M}$ of the type that appears in Proposition 
\ref{prop:findim}. 

We can write $\eta^C_{\lambdabold,M}$ as $d\log F^C_{\lambdabold,M}$, where $F^C_{\lambdabold,M}$ is
a  multivalued function (univalued if the exponents are integral) given by the product
\[
\prod_{\substack{i\in M\\1\le  \nu\le n}} (t_i-z_\nu)^{C(\alpha_{\bar i},\lambda_\nu)}
\!\!\!\!\prod_{i,j\in M, i\not=j} (t_i-t_j)^{-C (\alpha_{\bar i},\alpha_{\bar j})/2}
\prod_{\nu<\mu} (z_\nu-z_\mu)^{-C(\lambda_\nu,\lambda_\mu)}.
\]
We consider the first order differential operator $d^C:= d-\eta_{\lambdabold,M}^C$.
So $d^C$ is the ordinary exterior derivative for the multivalued polydifferentials after they get formally multiplied by the inverse of $F^C_{\lambdabold,M}$: $d^C= F^C_{\lambdabold,M}\,d\,(F^C_{\lambdabold,M})^{-1}$. In particular, a $d^C$-closed form is locally $F^C_{\lambdabold,M}$ times a $d$-closed form. This amounts to turning the trivial line bundle over $U_{n,M}$
into a local system $\LL_{\lambdabold,M}^C$: it is the local system for which $F^C_{\lambdabold,M}$  defines a flat (multivalued) section, in other words,
$\LL_{\lambdabold,M}^C$ is the structure sheaf of $\Ocal_{U_{n,M}}$ endowed with the connection for which $d^C$ is covariant derivation.

Let $\rho\in\hlie^*$ be, as usual in Lie theory, defined by the property that $\rho (\check{\alpha}_k)=1$  for $k=1,\dots ,r$ so that $C(\rho,\alpha_k)=q^C(\alpha_k)$ for every $k$.

\begin{lemma}\label{lemma:residuecalc}
The differential $\eta_{\lambdabold,M}^C$ has a logarithmic pole along each irreducible component of  $\Delta_{n,M}$ and we have
\begin{align*}
-\res_{\Delta (X)}\eta_{\lambdabold,M}^C &= q^C(\rho-\alpha_X)-q^C(\rho),\\
-\res_{\Delta_\infty(X)}\eta_{\lambdabold,M}^C &=q^C(\rho+\alpha_X)-q^C(\rho),\\
-\res_{\Delta_\nu(X)}\eta_{\lambdabold,M}^C &=q^C((\rho+\lambda_\nu)-\alpha_X)- q^C(\rho+\lambda_\nu),
\end{align*}
where $\alpha_X:=\sum_{x\in X} \alpha_{\bar x}$.
\end{lemma}
\begin{proof}
Since $\Delta (X)$ amounts to the confluence of the members of $X$, its generic point may be described in terms of $\Pbold^M\times U_n$ as the blow up of the diagonal obtained by setting all $t_x$, $x\in X$, equal to each other (followed by dividing out the action of the translation group of $\Abold$). The irreducible components of $D_{n,M}$ that
pass through this diagonal are defined by $t_x=t_y$, where $\{ x,y\}$ runs over the two-element subsets of $X$. The defining  formula for $\eta_{\lambdabold,M}^C$ shows that 
$\res_{(t_x=t_y)}\eta_{\lambdabold,M}^C=-\half C(\alpha_{\bar x},\alpha_{\bar y})$. We thus find that 
\begin{multline*}
-\res_{\Delta (X)}\eta_{\lambdabold,M}^C =\half \sum_{x\not=y}C(\alpha_{\bar x},\alpha_{\bar y})=\\
\half C(\sum_{x\in X}\alpha_{\bar x},\sum_{y\in X}\alpha_{\bar y})-
\half\sum_{x\in X}C(\alpha_{\bar x},\alpha_{\bar x})=\\
\half C(\alpha_X, \alpha_X)-\sum_{x\in X}C(\rho,\alpha_X)=
q^C(\alpha_X-\rho)-q^C(\rho).
\end{multline*}
In the case $\Delta_\infty(X)$, we also need to include the irreducible components of 
$D_{n,M}$ defined by $t_x=\infty$, $x\in X$. A straightforward computation shows that $-\res_{(t_x=\infty)}\eta_{\lambdabold,M}^C=C(\alpha_{\bar x},\alpha_{\bar x})$ and so we get  as additional term $\sum_{x\in X} C(\alpha_{\bar x},\alpha_{\bar x})=2 C(\rho,\alpha_X)$. This yields
$\half C(\alpha_X, \alpha_X)+\sum_{x\in X}C(\rho,\alpha_X)=
q^C(\alpha_X+\rho)-q^C(\rho)$, as asserted. Finally, for $\Delta_\nu(X)$, with $\nu=1,\dots ,n$, we need to subtract the residues for the divisors $t_x=z_\nu$, that is
$\sum_{x\in X} C(\alpha_{\bar x},\lambda_\nu)=C(\alpha_X,\lambda_\nu)$ and this gives the last value.
\end{proof}

The associated \emph{Aomoto complex} is the relative De Rham complex of $\Ocal_{U_n}$-modules $(p^+_*\Omega^\pt_{U^+_{n,M}/U_n}(\log \Delta_{n,M}), d^C)$.
Note that since logarithmic forms are $d$-closed, the relative differential $d^C$ is simply given by the wedge product with $-\xi_{\lambdabold,M}^C$.
Proposition  \ref{prop:findim} then tells us that:

\begin{lemma}
We have a natural isomorphism of $\Ocal_{U_n}$-modules
\[
\Ocal_{U_n}\otimes_\CC\tilde{V} (\lambdabold)_{0,\coprim}\cong
\Hcal^m(p^+_*\Omega^\pt_{U^+_{n,M}/U_n}(\log \Delta_{n,M}), d^C)\otimes_{\Scal (M_\pt)}\orient(M).
\]
In particular, $\Ocal_{U_n}\otimes_\CC V (\lambdabold)^\glie$ embeds in the latter.
\end{lemma}

The irreducible components of $\Delta_{n,M}$ along which the local system $\LL_{\lambdabold,M}^C$ has trivial monodromy are those on which the residue $\eta_{\lambdabold,M}^C$ is an integer. These play a special role in of a theorem of Esnault, Schechtman and Viehweg \cite{esv}, or rather our refinement (\cite{looijenga:kz}, \cite{looijenga:sl2}) which leads to a topological  interpretation  of the Aomoto cohomology.  For this purpose and for later uses, we  pause for a moment to discuss the various natural extensions of a rank on local system across a normal crossing divisor. 

\subsection*{Extensions across a normal crossing divisor}

Let $X$ be a complex manifold of dimension $n$, $D$ a simple normal crossing divisor on $X$ and  $\LL$  a rank one local system on $X-D$. If $D$ is irreducible (hence smooth), then there are two basic ways of extending $\LL$ to $X$ in the derived category:  extension by zero $R^\pt j_!\LL$ (where $j:X-D\subset X$ denotes the inclusion) and and the full direct image $R^\pt j_*\LL$. These are connected by a morphism $R^\pt j_!\LL(=j_!\LL)\to R^\pt j_*\LL$, 
 which on global cohomology this gives the  map $H^\pt (X,D;\LL)\to H^\pt (X-D;\LL)$. (Since  
$\LL$ is not defined on $D$, a word of explanation is in order: if $T$ is an open regular neighborhood of  
$D$ in $X$ relative to the Hausdorff topology so that its boundary $\p T$ lies in $X-D$, then 
$H^\pt(X,D;\LL)$ is understood as $H^m(X-T,\p T; \LL)$, or equivalently, 
as  cohomology with supports: $H^\pt_{\Phi}(X-D;\LL)$, where $\Phi$ is the collection of closed subsets of $X-D$ that
remain closed in $X-D$.)  If the monodromy of $\LL$ around $D$ is not the identity, then the two extensions coincide.  Shifted Verdier duality  converts this morphism of extensions into $R^\pt j_*\LL^\vee\leftarrow j_!\LL^\vee$.

Suppose  $D$ has two irreducible components $D'$ and $D''$ and we extended $\LL$ across the generic point of each of them so that we have an extension over $X-D'\cap D''$.  Then there is a natural extension over all of $X$ which is locally along a transversal slice of $D'\cap D''$  like an exterior product of two extensions as above over the complex unit disk. We can obtain it in stages, for instance,  by first doing the $D'$-extension over $X-D''$ and then the $D''$ extension over $X$; the opposite order yields the same result. More generally, if $D$ has several irreducible components, then an extension of $\LL$ over $X$ in the derived category is specified once we have done so at the generic points of $D$ and its formation is compatible with shifted Verdier duality.

So if we single out a collection of  irreducible components of $D$ along which $\LL$ has trivial monodromy and
denote its union $D^\sharp$, then we have specified an extension of $\LL$ over $X$: at the generic points of 
$D^\sharp$ we take the full direct image, and at the other generic points of $D$ we extend by zero. We denote that extension $\Ccal^\pt(\LL;D^\sharp)$. The cohomology of this extension is $H^\pt(X-D^\sharp,D-D^\sharp;\LL)$. 
Note however, that adding to $D^\sharp$ irreducible components of $D$ at which $\LL$ has nontrivial monodromy 
does no alter $\Ccal^\pt(\LL;D^\sharp)$ as an object in the derived category of constructible sheaves on $X$ and 
hence will  not affect the cohomology. The shifted Verdier dual of  $\Ccal^\pt(\LL;D^\sharp)$ is $\Ccal^\pt(\LL^\vee;D^\flat)$, where $D^\flat$ is the union of the irreducible components of $D$ at which $\LL$ has trivial monodromy not in $D^\flat$ (but we could equally well take the union of \emph{all} the irreducible components of $D$ not in $D^\sharp$). So if $X$ is of finite type, then we have a perfect pairing
\[
H^k(X, \Ccal^\pt(\LL;D^\sharp))\otimes H^{2n-k}_c(X,\Ccal^\pt(\LL^\vee;D^\flat))\to \CC.
\]

An important example is when $\LL\subset \Ocal_{X-D}$ is defined by a closed differential $\eta$. Given an  irreducible component  of $D$, then $\LL$ has trivial monodromy at it if and only if $\eta$ has there an integral residue. We
observed in \cite{looijenga:kz} and \cite{looijenga:sl2} that the twisted logarithmic De Rham complex  $(\Omega^\pt_{X}(\log D),d-\eta)$ represents $\Ccal^\pt (\LL_{\lambdabold,M}^C,D^{\ge 0})$, where $D^{\ge 0}$  is the union of the irreducible components of $D$ where $\eta$ has residue a nonnegative integer.

\subsection*{A Gau\ss-Manin connection}

The next proposition appears in \cite{looijenga:kz} and \cite{looijenga:sl2}, albeit that it is stated there in an absolute setting.  

\begin{proposition}\label{prop:topinterp}
Denote by $\Delta_{n,M}^{\ge 0}$ the union of the irreducible components of $\Delta_{n,M}$ where $\eta_{\lambdabold,M}^C$ has residue a nonnative integer. Then we have a natural identification of $\Ocal_{U_n}$-modules
\[
\Hcal^m(p^+_*\Omega^\pt_{U^+_{n,M}/U_n}(\log \Delta_{n,M}), d^C)\cong 
\Ocal_{U_n}\otimes R^m  p^+_*\Ccal^\pt (\LL_{\lambdabold,M}^C,\Delta_{n,M}^{\ge 0})
\] 
\end{proposition}
\begin{proof}
As noted above, the complex  $(\Omega^\pt_{U^+_{n,M}/U_n}(\log \Delta_{n,M}),d^C)$ represents the derived category object $\Ccal^\pt (\LL_{\lambdabold,M}^C,\Delta_{n,M}^{\ge 0})$.  The first lemma of section 2 of   \cite{esv} asserts that  the direct image
$R^q p^+_*\Omega^\pt_{U^+_{n,M}/U_n}(\log \Delta_{n,M})$ is zero unless  $q=0$. The proposition now follows by taking the $m$th  direct image on $U_n$.
\end{proof}

Note that the stalk at $\zz\in U_n$ of the sheaf  that appears  in the right hand side of the preceding 
proposition is equal to  the cohomology space
$H^m(U^+_{n,M}(\zz)-\Delta^{\ge 0}_{n,M},\Delta_{n,M}-\Delta^{\ge 0}_{n,M}; \LL_{\lambdabold,M}^C)$. 
The pair $(U^+_{n,M},\Delta_{n,M})$ is topologically locally trivial over $U_n$ and so 
$R^m  p^+_*\Ccal^\pt (\LL_{\lambdabold,M}^C,\Delta_{n,M}^{\ge 0})$ is a local system.  We conclude that the flat connection $d^C$ on $\Ocal_{U_{M,n}}$ induces one on  $\Ocal_{U_n}\otimes R^m  p^+_*\Ccal^\pt (\LL_{\lambdabold,M}^C,\Delta_{n,M}^{\ge 0})$ and (via  Proposition \ref{prop:topinterp}) one on  the $\Ocal_{U_n}$-module $\Hcal^m(p^+_*\Omega^\pt_{U^+_{n,M}/U_n}(\log \Delta_{n,M}), d^C)$. We will refer to this connection as  the \emph{Gau\ss-Manin connection}  and denote it by $\nabla_{GM}$.

\begin{corollary}\label{cor:main}
We  have a natural identification 
\[
\Ocal_{ U_n}\otimes_\CC\tilde{V} (\lambdabold)_{0,\coprim}\cong 
\Ocal_{U_n}\otimes R^m  p^+_*\Ccal^\pt (\LL_{\lambdabold,M}^C,\Delta_{n,M}^{\ge 0})\otimes_{\Scal (M_\pt)}\orient (M)
\] 
as $\Ocal_{ U_n}$-modules. (So this makes $\Ocal_{ U_n}\otimes V (\lambdabold)^\glie$ a trivial subbundle  of the right hand side.)
\end{corollary}

\subsection*{The KZ-connection is a Gau\ss -Manin connection}

Here is the main result:

\begin{theorem}\label{thm:gz=kz}
The embedding of  $\Ocal_{ U_n}\otimes V(\lambdabold)^\glie$ endowed with the KZ-connection 
$\nabla^C_{KZ}$ in  $\Hcal^m(p^+_*\Omega^\pt_{U^+_{n,M}/U_n}(\log \Delta_{n,M}), d^C)$ endowed with the  GM-connection $\nabla_{GM}$ is flat and hence induces an embedding of local systems $\KZ ^C(\lambdabold)\hookrightarrow R^m p^+_*\Ccal^\pt (\LL_{\lambdabold,M}^C,\Delta_{n,M}^{\ge 0})\otimes_{\Scal (M_\pt)}\orient (M)$.
\end{theorem}

Before we begin the proof, we show how  a bootstrap procedure turns this theorem into a more precise result.  

\begin{theorem}\label{thm:bettermain}
Let $\Delta_{n,M}^{> 0}$ denote the union of irreducible components of $\Delta_{n,M}$ along which the residue of  $\eta_{\lambda,M}^C$ is a positive integer. Then $\KZ ^{C}(\lambdabold)$ can be canonically identified with the isotypical part for the sign character of $\Scal (M_\pt)$ of the  image of  
\[
R^m p^+_*\Ccal^\pt(\LL_{\lambdabold,M}^{C},\Delta_{n,M}^{>0})\to R^m p^+_*\Ccal^\pt(\LL_{\lambdabold,M}^{C},\Delta_{n,M}^{\ge 0}).
\]
\end{theorem}

For clarity we note that the above map is at $\zz\in U_n$ the natural map of cohomology spaces with support
\[
H^m_{\Phi^{>0}}(U_{n,M}(\zz);\LL_{\lambdabold,M}^{C})\to 
H^m_{\Phi^{\ge 0}}(U_{n,M}(\zz);\LL_{\lambdabold,M}^{C}),
\]
where $\Phi^{>0}$ resp.\ $\Phi^{\ge 0}$ is the family of closed subsets of $U_{n,M}(\zz)$ which remain closed in 
$U^+_{n,M}(\zz)-\Delta_{n,M}^{>0}$ resp.\  $U^+_{n,M}(\zz)-\Delta_{n,M}^{\ge 0}$.

\begin{proof}
Theorem \ref{thm:gz=kz} yields an embedding
\[
\KZ ^{C}(\lambdabold)\hookrightarrow 
R^m p^+_*\Ccal^\pt(\LL_{\lambdabold,M}^{C},\Delta_{n,M}^{\ge 0})\otimes_{\Scal (M_\pt)}\orient (M).
\]
Next we note that the local system dual to $\KZ ^C(\lambdabold)$ is
$\KZ ^{-C}(\lambdabold')$, where the prime ${}'$ is the \emph{canonical involution} of $\hlie^*$, given as $-w_o$, where
$w_o$ is  the Weyl group element that maps the fundamental chamber  to its opposite. This
involution preserves the simple roots and (hence) the dominant weights.
With this notation, the involution applied to the identity $\lambda_1+\cdots +\lambda_n=m_1\alpha_1+\cdots +m_r\alpha_r$ yields $\lambda'_1+\cdots +\lambda'_n=m_1\alpha'_1+\cdots +\alpha'_r$ and so the data that served us for 
$\KZ ^C(\lambdabold)$ are for $\KZ ^{-C}(\lambdabold')$ given by $-C$ and the composite $\pi'$ of $\pi :M\to \{1,\dots ,r\}$ with the involution (also denoted by ${}'$) of $\{1,\dots ,r\}$  that is given by $\alpha'_k=\alpha_{k'}$.
We have thus defined $\eta^{-C}_{\lambdabold,\pi'}$. The definition then shows that
$\eta^{-C}_{\lambdabold',\pi' }=-\eta^C_{\lambdabold, M}$ so that  $\LL^{-C}_{\lambdabold',\pi'}$ may be identified with the dual of $\LL^{C}_{\lambdabold,M}$. Let us apply Corollary  \ref{cor:main} to the triple  $(-C,\lambdabold',\pi')$: we get a 
natural embedding of local systems
\[
\KZ ^{C}(\lambdabold)^\vee\hookrightarrow 
R^m p^+_*\Ccal^\pt(\LL_{\lambdabold,M}^{C},\Delta_{n,M}^{\le 0})\otimes_{\Scal (M_\pt)}\orient (M).
\]
By dualizing we obtain a surjection
\[
R^m p^+_*\Ccal^\pt(\LL_{\lambdabold,M}^{C},\Delta_{n,M}^{> 0})\otimes_{\Scal (M_\pt)}\orient (M)\twoheadrightarrow
\KZ ^{C}(\lambdabold)
\]
of local systems.  It remains to observe that the composite of the two relevant displays is the natural map. 
\end{proof}

We thus get a genuine topological characterization of the KZ-system.
For example, if $\eta^C_{\lambdabold,M}$ has no nonzero integral residues, then we find that $\KZ ^{C}(\lambdabold)_\zz$ can be identified with isotypical subspace of the sign character of the image of  $H^m_c(U_{n,M}(\zz),\LL_{\lambdabold,M}^C)\to H^m(U_{n,M}(\zz),\LL_{\lambdabold,M}^C)$. If in addition  $C$ is defined over $\RR$,  then $\LL_{\lambdabold,M}^C$ has  flat unitary metric (that gives $F^C_{\lambda,M}$ norm one) and there results hermitian (intersection) form on  this image, which is known to be nondegenerate. This puts on $\KZ ^{C}(\lambdabold)$ a (flat) nondegenerate hermitian form.

\begin{remark}
Theorem \ref{thm:gz=kz} tells us what the monodromy of the KZ system is like. If we fix  a base point  $\zz\in U_n$, then 
$\pi_1(U_n,\zz)$ is the colored braid group with $n$ strands. It acts on domain and range of the  linear map
$H^m_{\Phi^{>0}}(U_{n,M}(\zz);\LL_{\lambdabold,M}^{C})\to H^m_{\Phi^{\ge 0}}(U_{n,M}(\zz);\LL_{\lambdabold,M}^{C})$
in a manner that makes the map equivariant. It should be worthwhile to investigate such representations in their own right and perhaps make contact with the Kohno-Drinfeld approach via the representation theory of quantum groups. 
We further note that since the  KZ system embeds in a variation of complex mixed Hodge structure, it acquires a (flat) weight filtration. It should be interesting to  determine that filtration in terms of the KZ data. 
\end{remark}

\subsection*{When the Casimir element is defined over $\QQ$}
In case $C$ is defined over $\QQ$ in the sense that $q^C$ takes rational values on the roots, then 
$R^m \Ccal^\pt(\LL_{\lambdabold,M}^{C})$ is a eigen subsystem of a finite cyclic group acting on an ordinary  variation of mixed Hodge structure. To be precise, let $s$ be the smallest common denominator of these residues. Then the monodromy of $\LL_{\lambdabold,M}^{C}$ is the group $\bm{\mu}_s$ of $s$th roots of unity. It determines an unramified $\bm{\mu}_s$-cover  $\hat U_{n,M}\to  U_{n,M}$, so that the pull-back of $\LL_{\lambdabold,M}^C$ becomes trivial. This means that we may now $F^C_{\lambdabold,M}$ regard as a univalued (invertible) holomorphic function on 
$\hat U_{n,M}$.  Its normalization over $U^+_{n,M}$, $\hat U^+_{n,M}\to U^+_{n,M}$,  is a $\bm{\mu}_s$-cover of  $U^+_{n,M}$ that may have singularities, but these are quotient singularities and hence for our purposes of an innocent nature. The function $F^C_{\lambdabold,M}$ is meromorphic on it and the order of $F^C_{\lambdabold,M}$ along  an irreducible component of $\hat\Delta_{n,M}$ is given by $s$ times the residue of $\xi^C_{\lambdabold,M}$ along its image in $\Delta_{n,M}$ (indeed an integer).

If $\hat p^+: \hat U^+_{n,M}\to U_n$ denotes the projection,  then let $\hat\Delta_{n,M}^{\ge 0}, \hat\Delta_{n,M}^{>0},\dots $  have the obvious meaning. Then for $\zz\in U_n$, 
\[
(R^m_{\hat\Phi^{\ge 0}}\hat p_*\CC_{U^+_{n,M}})_\zz=
H^m(\hat U^+_{n,M}-\hat\Delta_{n,M}^{\ge 0}, \hat\Delta_{n,M}-\hat\Delta_{n,M}^{\ge 0};\CC)
\]
and  $R^m_{\hat\Phi^{\ge 0}}\hat p_*\CC_{U^+_{n,M}}$ comes with the structure of a variation of polarized mixed Hodge structure. There is now a finite group $\hat\Scal(M_\pt)$ acting on $\hat U^+_{n,M}$ which is an extension of  $\Scal (M_\pt)$ 
by the covering group $\mu_s$. It has a character $\hat\chi$ that is tautological on $\mu_s$ and lifts the sign character. 
On the level of stalks this yields the  identification of $\KZ ^{C}(\lambdabold)_\zz$ with the $\hat\chi$-isotypical space of
the image of $R^m_{\hat\Phi^{>0}}\hat p_*\CC_{U^+_{n,M}}\to R^m_{\hat\Phi^{\ge 0}}\hat p_*\CC_{U^+_{n,M}}$.

\subsection*{Proof of Theorem \ref{thm:gz=kz}}
We begin the proof  by recalling covariant derivation relative to the Gau\ss -Manin connection. Covariant derivation with respect to $z_\nu$ is exhibited on the form level by Lie derivation of a lift of this vector field to $U_{n,M}$. In order to ensure that logarithmicity is preserved we take a lift that depends on the argument: 

\begin{lemma}\label{lemma:gm}
Let $\zeta_{\Ibold}(\zz) =\zeta_{I^1}(z_1).\zeta_{I^2}(z_2).\cdots .
\zeta_{I^n}(z_n)$ be a basis element of  $\Bcal_{n,M}$ 
and let $\tilde\p_\nu:=\p_\nu+\sum_{i\in \{I^\nu\}}\frac{\p}{\p t_i}$ 
(a vector field on $\Pbold^{M}\times U_n$ that lifts the vector field  $\p_\nu$ to  $U_n$). Then the Lie derivative  $\Lcal^C_{\tilde\p_\nu}:=d^C\iota_{\tilde\p_\nu}+\iota_{\tilde\p_\nu}d^C$ maps $\zeta_{\Ibold}$ to  $-\xi_{\lambdabold,M}^C (\tilde\p_\nu).\omega$ and the latter lies in  $\CC[U_n]\otimes_\CC\Bcal_{n,M}$. This map is $\Scal_\pt$-equivariant and defines a connection on $\Bcal_{n,M}$ with logarithmic pole whose  form $A^C_{GM}$ lies in
\[
\sum_{\nu<\mu} \frac{d(z_\nu-z_\mu)}{z_\nu-z_\mu}\otimes_\CC \End_{\Abold^n} 
(\Bcal_{n,M}).
\]
We shall refer to this as the \emph{Gau\ss -Manin connection}.
\end{lemma}
\begin{proof}
We first notice that $\zeta_{\Ibold}$ is invariant  under the flow generated by $\tilde\p_\nu$ (which adds to the coordinates $(z_\nu, (t_i)_{i\in \{I^\nu\}})$ the same complex number), in other words,   $\Lcal_{\tilde\p_\nu}(\zeta_{\Ibold})=0$. Hence 
\begin{multline*}
\Lcal^C_{\tilde\p_\nu}(\zeta_{\Ibold})=d^C\iota_{\tilde\p_\nu}\zeta_{\Ibold} +\iota_{\tilde\p_\nu}d^C\zeta_{\Ibold}
=(d\iota_{\tilde\p_\nu}\zeta_{\Ibold} -\eta_{\lambdabold,M}^C\iota_{\tilde\p_\nu} \zeta_{\Ibold}) + \iota_{\tilde\p_\nu}(d\zeta_{\Ibold} - 
\eta_{\lambdabold,M}^C\zeta_{\Ibold})\\
= \Lcal_{\tilde\p_\nu}(\zeta_{\Ibold}) -\eta_{\lambdabold,M}^C\iota_{\tilde\p_\nu} \zeta_{\Ibold}- \iota_{\tilde\p_\nu} (\eta_{\lambdabold,M}^C\zeta_{\Ibold})=-\eta_{\lambdabold,M}^C (\tilde\p_\nu)\zeta_{\Ibold}.
\end{multline*}
For $i,j\in\{I^\nu\}$, the differentials $dt_i-dt_j$ and $dt_i-dz_\nu$ clearly vanish on 
$\tilde\p_\nu$ and so 
\begin{multline*}
-\eta_{\lambdabold,M}^C (\tilde\p_\nu)=
\sum_ {x\notin I^\nu} \frac{C(\alpha_{\bar x},\lambda_\nu)}{t_x-z_\nu}
-\sum_{i\in \{I^\nu\}}\sum_{\mu\not=\nu} \frac{C(\alpha_{\bar i},\lambda_\mu)}{t_i-z_\mu}\\
+\sum_{i\in \{I^\nu\}}\sum_{x\in \{I^\mu\}, \mu\not=\nu}\frac{C(\alpha_{\bar i},\alpha_{\bar x})}{t_i-t_x}
+\sum_{\mu\not=\nu}\frac{ C(\lambda_\nu,\lambda_\mu)}{z_\nu-z_\mu}.
\end{multline*}
Lemma's \ref{lemma:mixedshuffle} and \ref{lemma:mixedshuffle2} show that
multiplication of $\zeta_{\Ibold}$ by a  factor $(t_i-z_\mu)^{-1}$,  $(t_x-z_\nu)^{-1}$ or $(t_i-t_x)^{-1}$ ($i\in\{I^\nu\}$, $x\in \{I^\mu\}$, $\mu\not=\nu$) lands in $\frac{1}{z_\nu-z_\mu}\Bcal_n$. 
\end{proof}

 As is well-known (and easy to prove), $C$ has the form
\[
C=C_0+\sum_{\alpha} C_\alpha,
\]
with $C_0\in \hlie\otimes \hlie$ and $C_\alpha\in \glie_\alpha\otimes \glie_{-\alpha}$, where the sum is over all the roots. Here $C_0$ can be any symmetric tensor invariant under the Weyl group; it then  determines $C$. Since $C$ is symmetric,  $C_{-\alpha}$ is the transpose of $C_\alpha$.

We put 
\[
C_+:=\sum_{\alpha>0} C_\alpha\in \prod_{\alpha>0}\glie_\alpha\otimes \glie_{-\alpha}
 \quad \text{ and }\quad C_-:=\sum_{\alpha<0} C_\alpha\in \prod_{\alpha<0}\glie_\alpha\otimes \glie_{-\alpha},
\]
so that $C=C_++C_0+C_-$.
It is easy to check that $C_0$ acts semisimply in the tensor product of highest weight representations. In fact, for $\chi,\chi'\in\hlie^*$, $C_0$ acts on 
$V(\lambda)_{\chi}\otimes V(\lambda')_{\chi'}$ as multiplication by $C_0(\chi,\chi')$.
For the proof of Theorem \ref{thm:gz=kz}, we also need a better understanding of $C_+$.
The following lemma is essentially Lemma 7.6.3 of Schechtman and Varchenko \cite{sv},
and so we omit its proof.

\begin{lemma}\label{lemma:sv} 
Let $\tilde C_+:V(\lambda)\to V(\lambda)\otimes\glie$ be the linear map given by 
\[
\tilde C_+ (f_S 1_\lambda)=\sum_{\emptyset\not=T\le S} C(\alpha_{\ell (T)},\lambda-\alpha_{ S_{>\ell (T)}})f_{S-T}1_\lambda\otimes [f_T],
\]
where $\ell(T)$ denotes the last term of $T$ and the sum is taken over all nonempty subsequences of $T$ of $S$ and $S_{>\ell (T)}$ is the largest common tail of $S$  and $S-T$ (which of course may be empty). If $V$ is any representation of $\glie$, then the action of $C_+$ on $V (\lambda)\otimes V$ satisfies
\[
C_+ (f_S 1_\lambda \otimes v)=\tilde C_+ (f_S)(1_\lambda\otimes v).\qedhere
\]
\end{lemma}

\begin{proof}[Proof of Theorem \ref{thm:gz=kz}]
In view of the shape of the connections, it suffices to verify this in case $n=2$.
We begin working out  the computation in the proof of Lemma \ref{lemma:gm} in case $n=2$, $\nu=1$ (so that $\mu=1$).  We write $(z,w)$ for $(z_1,z_2)$,  $(\lambda,\mu)$ for $(\lambda_1,\lambda_2)$, $(I,J)$ for $(I^1,I^2)$ and $\zeta$ for $\zeta_I(z)\otimes\zeta_J(w)$. 

We have for $\tilde\p =\frac{\p}{\p z_\nu} +\sum_{i\in \{ I\}} \frac{\p}{\p t_i}$:
\begin{multline*}
-\eta_{\lambdabold,M}^C (\tilde\p)=
\sum_ {j\in \{ J\}} \frac{C(\alpha_{\bar{j}},\lambda)}{t_{j}-z}
-\sum_{i\in \{ I\}} \frac{C(\alpha_{\bar i},\mu)}{t_i-w}
+\sum_{\substack{i\in \{ I\}\\ j\in \{ J\}}}\frac{C(\alpha_{\bar i},\alpha_{\bar{j}})}{t_i-t_{j}}
+\frac{C(\lambda,\mu)}{z-w}
\end{multline*}
and so 
\begin{multline*}
(z-w)\eta_{\lambdabold,M}^C (\tilde\p)\zeta=
-\sum_ {j\in \{ J\}} C(\alpha_{\bar{j}},\lambda)\frac{w-z}{t_{j}-z}\zeta_I(z)\otimes\zeta_J(w)\\
-\sum_{i\in \{ I\}} C(\alpha_{\bar i},\mu)\frac{z-w}{t_i-w}\zeta_I(z)\otimes\zeta_J(w)\\
+\sum_{i\in \{ I\},j\in \{ J\}} C(\alpha_{\bar i},\alpha_{\bar{j}})\frac{z-w}{t_i-t_{j}}\zeta_I(z)\otimes\zeta_J(w)
- C(\lambda,\mu)\zeta.
\end{multline*}
We develop these terms with the help of Lemmas \ref{lemma:mixedshuffle} and \ref{lemma:mixedshuffle2} and get
\begin{multline*}
\sum_{\substack{I=I''iI'\\ J=J''jJ'}} C(\alpha_{\bar i},\alpha_{\bar{j}})\cdot
\Big(\zeta-\sum_{I''=I_2I_1} (-1)^{|I_2|} \omega_{I_1}(z)
\otimes\omega_{I_2^*i}\omega_{I''i}\frac{dt_i}{t_i-t_{j}}\zeta_J(w)\\
-\sum_{J'=J_2J_1} (-1)^{|J_2|}
\omega_{J_2^*j}\omega_{J''j}\frac{dt_{j}}{t_{j}-t_i}\zeta_I(z)\otimes\zeta_{J_1}(w)\Big)\\
-\sum_{I=I''iI'} C(\alpha_{\bar i},\mu)
\Big(\zeta-\sum_{I'=I_2I_1} (-1)^{|I_2|}\zeta_{I_1}(z)
\otimes\omega_{I_2^*i)}\zeta_{I''i}(w)\zeta_J(w)\Big)\\
-\sum_ {J=J''jJ'} C(\alpha_{\bar{j}},\lambda)\Big(\zeta
-\sum_{J'=J_2J_1} (-1)^{|J_2|}
\omega_{J_2^*j)}\omega_{J''j}(z)\zeta_I(z)
\otimes\zeta_{J_1}(w)\Big)+\\+C(\lambda,\mu),
\end{multline*}
which after collecting terms becomes
\begin{multline*}
C\big(\lambda-\sum_{i\in \{ I\}}\alpha_{\bar i}, \mu-\sum_{j\in \{ J\}} \alpha_{\bar{j}}\big)\zeta
\\
+\sum_{I=I''iI_2I_1}(-1)^{|I_2|} \zeta_{I_1}(z)\otimes
\omega_{I_2^*i}\omega_{I''i}
\Big(\sum_{j\in \{ J\}} \frac{C(\alpha_{\bar i},\alpha_{\bar{j}})}{t_i-t_{j}}-
\frac{C(\alpha_{\bar i},\mu)}{t_i-w}\Big)dt_i\zeta_J(w)\\
+\!\!\!\sum_{J=J''jJ_2J_1}\!\!(-1)^{|J_2|} \omega_{J_2^*j}\omega_{J''j}\cdot
\Big(\sum_{i\in N} \frac{C(\alpha_{\bar i},\alpha_{\bar{j}})}{t_{j}-t_i}
-\frac{C(\lambda,\alpha_{\bar{j}})}{t_{j}-z}\Big)dt_{j}\zeta_I(z)\otimes \zeta_{J_1}(w).
\end{multline*}
Since
\[
\Big(\sum_{j\in \{ J\}} \frac{C(\alpha_{\bar i},\alpha_{\bar{j}})}{t_i-t_{j}}-
\frac{C(\alpha_{\bar i},\mu)}{t_i-w}\Big)dt_i\zeta_J(w)=-\half C(\alpha_{\bar{i}},\alpha_{\bar{i}}) \Phi_i(\zeta_J(w)),
\]
we may also write the previous expression as $\G_0(\zeta)+\G_+(\zeta_I(z)) \zeta_J(w)+
\G_-(\zeta_J(w))\zeta_I(z)$ with
$\G_0(\zeta)=C_0(\lambda-\sum_{i\in I}\alpha_{\bar i}, \mu-\sum_{j\in J} \alpha_{\bar{j}})\zeta$
and where $\G_+(\zeta_I(z))$ is the operator from $\Bcal $ to $\Bcal \otimes\Bcal .1_\mu$ defined by
\[
\G_+(\zeta_I(z))=
\sum_{I=I''iI_2I_1} (-1)^{|I_2|} \half C(\alpha_{\bar{i}},\alpha_{\bar{i}}) \zeta_{I_1}(z)\otimes\omega_{I''i}\omega_{I_2^*i}\Phi_{i},
\]
and $\G_-$ is its transpose. It follows from Corollary \ref{cor:funny} of the appendix that $\G_+$ defines a linear
map $\G :V(\lambda)\to V(\lambda)\otimes_\CC\glie$ that coincides with the map $\tilde C_+$ that appears in Lemma \ref{lemma:sv}. Since $\G_-$ resp.\ $C_-$ is the transpose of $\G_+$ resp.\ $C_+$, we
conclude that identity  we were after indeed holds: $A^C_{GM}= C_{1,2}\otimes d(z-w)/(z-w)=A^C_{KZ}$.
\end{proof}

\section{The WZW-system}

Our discussion of the case when $C$ is defined over $\QQ$ covers one of particular interest, namely  the one for which is defined the WZW-subsystem of a given level, where it is assumed that $\glie$ is simple and finite dimensional. We recall its definition. Let $\tilde\alpha\in\hlie^*$ be the highest root relative to the root basis $(\alpha_1,\dots ,\alpha_r)$ and let $\tilde{\alpha}^\vee$ be the associated coroot. We fix a generator $e$ of the (one dimensional) root space $\glie_{\tilde\alpha}$ and define an $\Ocal_{U_n}$-linear endomorphism $\Ecal$ of  $\Ocal_{U_n}\otimes V(\lambdabold)$ by 
\[
\Ecal(\zz)=\sum_{\nu=1}^n 1\otimes \cdots \otimes 1\otimes (z_\nu e)\otimes 1\otimes\cdots \otimes 1,
\]
where $z_\nu e$ is acting on $V(\lambda_\nu)$. So if for all $\nu$, $z_\nu\not=0$ (a property we can arrange for by doing a translation in $\Abold$), and we let $\glie$ act on $V(\lambda_\nu)$ by  
modifying the given action in terms of the scalar $z_\nu$: $e(z_\nu)_k:=z_\nu e_k$ resp.\  $f(z_\nu)_k:=z_\nu^{-1}f_k$, then $\Ecal(\zz)$ acts on  $V(\lambdabold)$ as $e$. 

\begin{lemma}\label{lemma:Echar}
Let $I$  be a  sequence in $\Ical $ representing the highest root. Then for a suitable choice of $e$, we have that for every
$\zeta\in \Ocal_{U_n}\otimes V(\lambdabold)$, $(\Ecal (\zeta))_{\Ical-\{ I\}}=\res_{t_I=\infty}t_I\res_I$. 
\end{lemma}
\begin{proof}
Let $\zeta\in V(\lambdabold)$ be of the form $\zeta (\zbold)=
\zeta_1(z_1)\cdots \zeta_n(z_n)$ with  $\zeta_\nu\in V(\lambda_\nu)$. Then 
$e (\zeta_\nu(z_\nu))_{\Ical -\{ I\}}= \res_{(t_I=z_\nu)}\res_I\zeta_\nu(z_\nu)$ by Corollary
\ref{cor:diagonal2} and so
$z_\nu e (\zeta_\nu(z_\nu))_{\Ical -\{ I\}}=\res_{(t_I=z_\nu)}t_I\res_I\zeta_\nu(z_\nu)$. It follows that
$\Ecal(\zeta (\zbold))_{\Ical -\{ I\}}=\sum_{\nu=1}^n \res_{(t_I=z_\nu)}t_I\res_I\zeta(\zbold)=-\res_{(t_I=\infty)}t_I\res_I\zeta(\zbold)$.
\end{proof}

Let $\ell$ be a fixed positive integer.  We say that a representation $V$ of $\glie$ \emph{is of level $\le \ell$} if  $e^{1+\ell}$ is the zero endomorphism in $V$. For $V=V(\lambda)$, with $\lambda$ dominant integral, this means that  $\lambda (\tilde{\alpha}^\vee)\le \ell$.  In what follows we assume that our $V(\lambda_1),\dots ,V(\lambda_n)$ are all of level $\le \ell$. According to Corollary \ref{cor:diagonal2} this amounts to the condition that for any $(1+\ell)$-tuple of sequences $(I_0,\dots ,I_{\ell})$ in $\Ical $ representing the highest root,  we have 
\[
\res_{(t_{I_0}=z_\nu)}\cdots \res_{(t_{I_{\ell}}=z_\nu)}\res_{I_0}
\res_{I_1}\cdots \res_{I_{\ell}}\zeta=0,\quad (\nu=1,\dots ,n).
\]  

The \emph{Verlinde space of level $\ell$}  is defined in a setting which involves a punctured compact Riemann surface as its `continuous input' so that over the moduli space of such punctured  Riemann surfaces these spaces make up a vector bundle,  the so-called \emph{WZW-bundle of level $\ell$}. (When this bundle is pulled back to a certain $\CC^\times$-bundle over that  moduli space, it acquires a natural flat connection.)  In case of the Riemann sphere, the sheaf of sections of this bundle (or of its dual, depending on convention) may be obtained as a subbundle $\Wcal_\ell (\lambdabold)$ of $\Ocal_{U_n}\otimes  V(\lambdabold)^\glie$:
\[
\Wcal_\ell (\lambdabold):=\ker (\Ecal^{1+\ell}| \Ocal_{U_n}\otimes V(\lambdabold)^\glie).
\]
So its fiber over $\zbold$ yields the vectors in $V(\lambdabold)^\glie$ that generate a representation of level $\le\ell$ relative to  the modified $\glie$-representations on the factors. 

We recall that the length of the highest root $\tilde\alpha$ is one less than the  \emph{Coxeter number} $h$ of $\glie$.  

\begin{corollary}\label{cor:wzwvanishing}
An element of  $\zeta\in \Ocal_{U_n}\otimes V(\lambdabold)^\glie$ lies in $\Wcal_\ell (\lambdabold)$
if and only if for any $(1+\ell)$-tuple of sequences $(I_0,\dots ,I_{\ell})$ in $\Ical $ with each 
member  of length $h-1$, $\res_{I_0}\res_{I_1}\cdots \res_{I_{\ell}}\zeta$ vanishes on the diagonal locus 
defined by $\cup_k \{I_k\}\subset \Ical$.
\end{corollary}
\begin{proof}
If $I$ is a finite sequence in $\Ical$ such that $\res_I$ is nonzero on $V(\lambda_\nu)$, then
$I$ has length $\le h-1$ and in case of equality, $I$ represents the highest root $\tilde\alpha$.
Lemma \ref{lemma:Echar} tells us that  $\zeta\in \Ocal_{U_n}\otimes V(\lambdabold)^\glie$ lies in $\Wcal_\ell (\lambdabold)$ if and only if for any $(1+\ell)$-tuple of sequences $(I_1,\dots ,I_{1+\ell})$ in $\Ical $ representing the highest root,  then 
\[
\res_{(t_{\{I_0\}}=\infty)}\!\cdots\! \res_{(t_{\{I_{\ell}\}}=\infty)}\!\big(t_{\{I_0\}}t_{\{I_1\}}\cdots t_{\{I_{\ell}\}}\!\res_{I_0}\!
\res_{I_1}\cdots \res_{I_{\ell}} \zeta\big)\! =\!0.
\]  
So if we put  $\omega:=\res_{I_0}\res_{I_1}\cdots \res_{I_{\ell}} \zeta$ (a  polydifferential on the diagonal  with coordinates $t_{\{I_0\}},\dots ,t_{\{I_\ell\}}, \{ t_j\}_{j\in M-\{ I\}}$), then it remains to see that the  property is equivalent to the vanishing of $\omega$ at $t_{\{ I_0\}}=\cdots =t_{\{ I_{\ell }\}}$. Since $\zeta$ is $\aut (\Pbold)$-invariant, so is  $\omega$, and hence it suffices to prove that the latter vanishes on $t_{\{ I_0\}}=\cdots =t_{\{ I_{\ell }\}}=\infty$.

Recall that by  Theorem \ref{thm:diagonaltensor} the polar loci of $\omega$ involving
the coordinates $t_{\{I_0\}},\dots ,t_{\{I_{\ell}\}}$, are of the type  $(t_{\{I_k\}}=z_\nu)$ only. 
Now let us put for $k=0,\dots , \ell$, $u_k:=t_{\{ I_k\}}^{-1}$. So  $\omega$ is regular in the generic point defined by $u_1=\cdots =u_{1+\ell}=0$. The vanishing condition amounts to: 
\[
\res_{(u_0=0)}\cdots \res_{(u_{\ell}=0)} \frac{\omega}{u_0u_1\cdots u_{\ell }}=0.
\]
This is equivalent to: $\omega$  vanishes on the locus  $u_0=\cdots =u_{\ell}=0$.
\end{proof}

In case  $\glie=\slin (2)$, we have $r=1$ and the highest root is the unique simple root (so that $h-1=1$). Corollary 
\ref{cor:wzwvanishing} then says that $\zeta\in \Ocal_{U_n}\otimes V(\lambdabold)^\glie$ lies in $\Wcal_\ell (\lambdabold)$ if and only if $\zeta$ vanishes on any diagonal defined by an $(1+\ell)$-element subset of $\Ical$. This is  due to Ramadas \cite{ramadas}, who proved this  in an entirely different manner. The proof given here is closer in spirit to ours in \cite{looijenga:sl2}.

Beilinson and Feigin have shown that $\Wcal_\ell (\lambdabold)$ is locally free and flat for the KZ connection $\nabla^{C_\ell}_{KZ}$, where $C_\ell$ is  characterized by the fact that $q^{C_\ell}(\tilde \alpha)=(\check{h}+\ell)^{-1}$, where $\check{h}:=1+\rho (\tilde{\alpha}^\vee)$ is known as the \emph{dual Coxeter number} of $\glie$ (strictly speaking, they prove the dual statement). In particular, $C_\ell$ is defined over $\QQ$.  A long standing conjecture in the physicists's community is the existence of a flat unitary metric on this bundle.

\begin{conjecture}\label{conj:main}
The subbundle $\Wcal_\ell (\lambdabold)$  maps to the square integrable forms, or what amounts to the same,  lands in 
sign isotypical part of the direct image of the relatively dualizing sheaf, i.e., in $\hat p^+_*\omega_{\hat U^+_{n,M}/U_n})$.  In particular,  it is of pure bidegree $(m,0)$ and the WZW system has a flat unitary structure.
\end{conjecture}

For $\glie=\slin (2)$ this has been proved by Ramadas \cite{ramadas}, who
derives it from  the above vanishing property on codimension $\ell$ diagonals (see also \cite{{looijenga:sl2}}). 

\section{Appendix: an operator formula}
We take up the situation of Section \ref{sect:polydiff}, but will assume  that $c_{i,j}=c_{j,i}$.
Recall that if $I$ is a sequence, then $I^*$ denotes that sequence in reverse order.

\begin{lemma}\label{lemma:funny}
Assume  that $c_{i,j}=c_{j,i}$ and let $\G : \Bcal\to \Bcal\otimes_\CC\End (\Bcal )$ be the linear map defined by 
\[
\G (\zeta_I):=\sum_{I=LiKJ} (-1)^{|K|}\zeta_J\otimes
\omega_{Li}\omega_{K^*i}\Phi_i
\]
(so $\G (1)=0$ and $\G(\zeta_i)=1\otimes\Phi_i$). Then for any finite sequence $I$ in $\Ical$ and $x\in\Ical$, we have 
\[
\G(\Phi_x(\zeta_I))=(\Phi_x\otimes 1 +1\otimes\ad \Phi_x)\G(\zeta_I)
+(p_x-c_{x,I})\zeta_I\otimes \Phi_x.
\]
\end{lemma}

\begin{corollary}\label{cor:funny}
In this situation we have 
\[
\G (\Phi_I (1))=
\sum_{\emptyset\not=K\le I} (p_{\ell(K)}-c_{\ell (K),I_{>\ell (K)}})\Phi_{I-K}(1)\otimes [\Phi_K],
\]
where the sum is over all nonempty subsequences of $I$, where $\ell(K)$ denotes
the last term of $K$ and $I_{>\ell (K)}$ is the largest common tail of $I$ and $I-K$.
\end{corollary}
\begin{proof}
We have that $\Phi_I(1)$ is a linear combination of the $\zeta_{\sigma (I)}$, where $\sigma$ runs over the permutations of $I$. So in the preceding lemma we may replace $\zeta_I$ by $\Phi_{I}(1)$. Then the  claimed identity  follows with induction.
\end{proof}

\begin{proof}[Proof of Lemma \ref{lemma:funny}]
The proof is a  mixture of algebra and bookkeeping and not entirely straightforward.

We derived  in Lemma \ref{lemma:betaop} the identity
\[
\Phi_x(\zeta_I)
=\sum_{I=I''I'} (p_x-c_{x,I'})\zeta_{I''xI'}.
\]
So each term that appears in $\G (\Phi_x\zeta_I)$ corresponds to a way of
writing  $I''xI'$ as $\tilde Li\tilde K\tilde J$. We can also express this differently by writing 
$I$ as $LiKJ$, and then insert $x$ in resp.\ $J$, $K$, $L$, or write $I=LKJ$ and 
take  $LxKJ$ (this is when $i=x$). We thus get
\begin{align*}
\G (\Phi_x&(\zeta_I))=\\
&\sum_{I=LiKJ, J=J''J'} (-1)^{|K|}(p_x-c_{x,J'})\zeta_{J''xJ'}\otimes \omega_{Li}\omega_{K^*i}\Phi_i\, +\tag{I}\\
&\sum_{I=LiKJ, K=K''K'} -(-1)^{|K|}(p_x-c_{x,K'J})\zeta_J\otimes \omega_{Li}\omega_{(K')^*x(K'')^*i}\Phi_i\, +\tag{II}\\
&\sum_{I=LiKJ, L=L''L'} (-1)^{|K|}(p_x-c_{x,L'iKJ})\zeta_J\otimes \omega_{L''xL'i}\omega_{K^*i}\Phi_i\, +\tag{III}\\
&\sum_{I=LKJ} (-1)^{|K|}(p_x-c_{x,KJ})\zeta_J\otimes \omega_{Lx}\omega_{K^*x}\Phi_x.\tag{IV}
\end{align*}
Denoting the subsums appearing above by their roman tags, then we observe that $\text{(I)}=(\Phi_x \otimes 1) \G (\zeta_I)$. We rework $\text{(IV)}$ by writing it first as 
\[
\text{(IV)}=\sum_{I=MJ}\zeta_J\otimes\Big(\sum_{M=LK} (-1)^{|K|}(p_x-c_{x,J}-c_{x,K})\omega_{Lx}\omega_{K^*x}\Big)\Phi_x
\]
and then continue with the inner sum. The expression $\omega_{Lx}\omega_{K^*x}$ can be written as shuffle product
$\sum_{S\in L\star K^*} \omega_{Sx}$. In case $S$ is empty, this reduces to just $\omega_x$, but otherwise such a shuffle appears twice: if $S$ ends with $i$ and $M$ is written $M''iM'$, then either $K=M'$ or $K=iM'$. These terms appear with coefficients $(-1)^{|M'|}(p_x-c_{x,J}-c_{x,M'})$ resp.\ $(-1)^{|M'|+1}(p_x-c_{x,J}-c_{x,M'}-c_{x,i})$ and 
so add up to give $(-1)^{|M'|}c_{x,i}$.
We conclude that (after substituting $L$ for $M''$ and $K$ for $M')$: 
\[
\text{(IV)}=(p_x-c_{x,I})\zeta_I\otimes \omega_x\Phi_x+
\sum_{I=LiKJ}(-1)^{|K|}c_{x,i}\zeta_J\otimes\omega_{Lix}\omega_{K^*ix}\Phi_x.
\]

Next we compute 
\[
(1\otimes \ad (\Phi_x))\G(\zeta_I)=
\sum_{I=LiKJ} (-1)^{|K|}\zeta_J\otimes [\Phi_x,\omega_{Li}\omega_{K^*i}\Phi_i]
\]

We work out the expression  $[\Phi_x,\omega_{Li}\omega_{K^*i}\Phi_i]$ on the right of the tensor symbol;
it is the sum 
\[
[\Phi_x,\omega_{Li}]\omega_{K^*i}\Phi_i+\omega_{Li} [\Phi_x,\omega_{K^*i}]\Phi_i+
\omega_{Li}\omega_{K^*i} [\Phi_x,\Phi_i]
\]
 and for these terms we have according to Lemma
\ref{lemma:betaop},
\begin{align*}
 [\Phi_x,\omega_{Li}]\omega_{K^*i}\Phi_i&=\sum_{L=L''L'} -c_{x,L'} \omega_{L''xL'i}\omega_{K^*i}\Phi_i,\\
\omega_{Li} [\Phi_x,\omega_{K^*i}]\Phi_i&=\sum_{K=K''K'} -c_{x,K''}
\omega_{Li}\omega_{(K')^*x(K'')^*i}\Phi_i.
\end{align*}
According to Lemma \ref{lemma:simplebracket}, 
$[\Phi_x,\Phi_i]=-c_{x,i}\omega_{(x,i)}\Phi_i +c_{i,x}\omega_{(i,x)}\Phi_x$ and hence we  find
\begin{multline*}
\omega_{Li}\omega_{K^*i} [\Phi_x,\Phi_i]=-c_{x,i}\omega_{Li}\omega_{K^*i}\omega_{(x,i)}\Phi_i
+c_{i,x}\omega_{Li}\omega_{K^*i}\omega_{(i,x)}\Phi_x=\\
-c_{x,i}\omega_{Li}\omega_{K^*i}\omega_{(x,i)}\Phi_i +c_{i,x}\omega_{Lix}\omega_{K^*ix}\Phi_x.
\end{multline*}
It follows that
\begin{align*}
(1\otimes \ad (&\Phi_x))\G(\zeta_I)=\\
&\sum_{I=LiKJ,L=L''L'} -(-1)^{|K|}c_{x,L'}\zeta_J\otimes\omega_{L''xL'i}\omega_{K^*i}\Phi_i\, +\tag{V}\\
&\sum_{I=LiKJ, K=K''K'} -(-1)^{|K|}c_{x,K''}\zeta_J\otimes\omega_{Li}
\omega_{(K')^*x(K'')^*i}\Phi_i\,+\tag{VI}\\
&\sum_{I=LiKJ} -(-1)^{|K|}c_{x,i}\zeta_J\otimes \omega_{Li}\omega_{K^*i}\omega_{(x,i)}\Phi_i\,+\tag{VII} \\
&\sum_{I=LiKJ} (-1)^{|K|}c_{i,x}\zeta_J\otimes\omega_{Lix}\omega_{K^*ix}\Phi_x.\tag{VIII}
\end{align*}
Adhering to our convention of identifying subsums by the corresponding roman tags,  we see that
\begin{multline*}
\text{(II)}-\text{(VI)}\!\! =\!\!\sum_{I=LiKJ}\!\! - (-1)^{|K|}(p_x-c_{x,KJ})\zeta_J\otimes\omega_{Li}\sum_{K=K''K'}\omega_{(K')^*x(K'')^*i}\Phi_i\\
=\sum_{I=LiKJ}  -(-1)^{|K|}(p_x-c_{x,KJ})\zeta_J\otimes\omega_{Li}\omega_{(x,i)}\omega_{K^*i}\Phi_i 
\end{multline*}
and similarly
\begin{multline*}
\text{(III)}-\text{(V)}-\text{(VII)}\! =\!\!\sum_{I=LiKJ}\!\!\! (-1)^{|K|}(p_x-c_{x,KJ})\zeta_J\otimes\sum_{L=L''L'}\omega_{L''xL'i}\omega_{K^*i}\Phi_i\\
=\sum_{I=LiKJ} (-1)^{|K|}(p_x-c_{x,KJ})\zeta_J\otimes \omega_{Li}\omega_{(x,i)}\omega_{K^*i}\Phi_i,
\end{multline*}
so that $\text{(II)}+\text{(III)}-\text{(V)}-\text{(VI)}-\text{(VII)}=0$. 

Since $c_{x,i}=c_{i,x}$, we see that $\text{(IV)}-\text{(VIII)}=(p_x-c_{x,I})\zeta_I\otimes \Phi_x$.
The Lemma follows.
\end{proof}


\begin{thebibliography}{1}

\bibitem{bd}
A.~Beilinson, V.G.~Drinfeld: 
\textit{Affine Kac-Moody algebras and polydifferentials,}  
Internat.\ Math.\ Res.\ Notices 1, 1--11  (1994).

\bibitem{drinfeld}
V. G.~Drinfeld:
\textit{Quasi-Hopf algebras and Knizhnik-Zamolodchikov equations.} 
Problems of modern quantum field theory (Alushta, 1989), 1Ð13, 
Res. Rep. Phys., Springer, Berlin, 1989. 

\bibitem{esv}
H.~Esnault, V.~Schechtman, E.~ Viehweg:
\textit{Cohomology of local sytems on the complement of hyperplanes},
Invent.\ Math.\ 109, 557--661 (1992)
\textit{Erratum}, Invent.\ Math.\ 112, 447 (1993).

\bibitem{kac}
V.~ Kac:
\textit{Infinite dimensional Lie algebras}, 3rd edition, Cambridge University Press,
Cambridge (UK), 1990.

\bibitem{looijenga:kz}
E.~Looijenga:
\textit{Arrangements, KZ systems and Lie algebra homology,} Singularity theory (Liverpool, 1996), xvi, 109--130, London Math. Soc. Lecture Note Ser., 263, Cambridge Univ. Press, Cambridge (UK), 1999.

\bibitem{looijenga:sl2}
E.~Looijenga:
\textit{Unitarity of $\SL (2)$-conformal blocks in genus zero,} 
J.\ Geom.\ Phys.\ 59, 654--662 (2009).

\bibitem{ramadas}
T.R.~Ramadas: 
\textit{The ``Harder-Narasimhan trace'' and unitarity of the KZ/Hitchin connection: genus 0,}  
Ann.\ of Math.\ (2)  169 1--39 (2009).

\bibitem{sv} 
V.V.~Schechtman, A.N.~Varchenko:
\textit{Arrangements of hyperplanes and Lie algebra homology},
Invent.\ Math.\ 106, 139--194 (1991). 

\bibitem{sf}
A.V.~Stoyanovski, B.L.~Feigin: 
\textit{Realization of a modular functor in the space of differentials, and geometric approximation of the manifold of moduli of $G$-bundles,}  Funktsional.\ Anal.\ i Prilozhen.\  28  (1994), 42--65, 95; translation in  Funct.\ Anal.\ Appl.\ 28 257--275 (1995).


\end{thebibliography}
\end{document}